\documentclass[11pt]{amsart}

\usepackage{tikz}
\usetikzlibrary {positioning}

\usepackage[square,compress,comma, numbers,sort]{natbib}
\usepackage[colorlinks=true, citecolor=red, linkcolor=blue]{hyperref}
\usepackage{amsfonts,mathtools}

\allowdisplaybreaks[4]

\usepackage{amsmath}
\usepackage{bbm}
\usepackage{amssymb}
\usepackage{mathtools}
\usepackage{leftindex}

\usepackage{mathrsfs}
\usepackage[scr=rsfs,cal=boondox]{mathalfa}

\usepackage{color}

\definecolor{c20}{rgb}{0.,0.7,0.}
\definecolor{c30}{rgb}{0.,0.,1.}
\definecolor{c40}{rgb}{1,0.1,0.7}
\definecolor{c50}{rgb}{1,0,0}
\definecolor{c60}{rgb}{1,0.9,0.1}

\newcommand{\abs}[1]{\left\lvert #1 \right\rvert}

\newcommand{\sprod}[1]{\langle#1\rangle}

\newcommand{\E}[1]{\mathbb{E}\left\{ #1\right\}}

\newcommand{\R}{\mathbb{R}}

\newcommand{\N}{\mathbb{N}}

\newcommand{\ldot}{,\ldots,}

\newcommand{\BQN}{\begin{eqnarray}}
\newcommand{\EQN}{\end{eqnarray}}
\newcommand{\BQNY}{\begin{eqnarray*}}
\newcommand{\EQNY}{\end{eqnarray*}}

\newcommand{\BS}{\begin{sat}}
\newcommand{\ES}{\end{sat}}
\newcommand{\BT}{\begin{theo}}
\newcommand{\ET}{\end{theo}}
\newcommand{\BK}{\begin{korr}}
\newcommand{\EK}{\end{korr}}

\newcommand{\BD}{\begin{de}}
\newcommand{\ED}{\end{de}}

\newcommand{\BIT}{\begin{itemize}}
\newcommand{\EIT}{\end{itemize}}
\newcommand{\BDI}{\begin{description}}
\newcommand{\EDI}{\end{description}}

\newcommand{\BRM}{\begin{remarks}}
\newcommand{\ERM}{\end{remarks}}

\newcommand{\BEL}{\begin{lem}}
\newcommand{\EEL}{\end{lem}}

\newtheorem{theo}{Theorem}[section]
\newtheorem{sat}[theo]{Proposition}
\newtheorem{de}[theo]{Definition}
\newtheorem{lem}[theo]{Lemma}

\newtheorem{example}[theo]{Example}
\newtheorem{korr}[theo]{Corollary}
\newtheorem{remark}[theo]{Remark}
\newtheorem{remarks}[theo]{Remarks}

\newcommand{\nelem}[1]{{Lemma \ref{#1}}}

\newcommand{\netheo}[1]{{Theorem \ref{#1}}}

\newcommand{\prooflem}[1]{\textsc{\bf Proof of Lemma} \ref{#1}:}

\newcommand{\COM}[1]{}

\def\td{\text{\rm d}}

\newcommand{\QED}{\hfill $\Box$}

\newcommand{\kb}[1]{\boldsymbol{#1}}
\newcommand{\vk}[1]{\kb{#1}}

\title[L\'{e}vy-driven queuing networks in 
multi-scale light and heavy traffic]{L\'{e}vy-driven queuing networks in \\
multi-scale light and heavy traffic}

\author[K. D\c{e}bicki]{Krzysztof D\c{e}bicki}
\address{Krzysztof D\c{e}bicki, Mathematical Institute\\
        Wroc\l{}aw University\\
        pl.\ Grunwaldzki 2/4, 50-384 Wroc\l{}aw, Poland}
\email{krzysztof.debicki@math.uni.wroc.pl}

\author[N. Kriukov]{Nikolai Kriukov}
\address{Nikolai Kriukov, Korteweg-de Vries Institute\\
	University of Amsterdam\\
	P.O. Box  94248, 1090 GE  Amsterdam, The Netherlands
}
\email{n.kriukov@uva.nl}

\author[M. Mandjes]{Michel Mandjes}
\address{Michel Mandjes, Mathematical Institute\\
        Leiden University\\
        P.O. Box 9512, 2300 RA Leiden, The Netherlands}
\email{m.r.h.mandjes@math.leidenuniv.nl}

\begin{document}

\maketitle

\begin{abstract}
   We study a queueing network with a strictly upper-triangular routing matrix, 
   where each column contains at most one non-negative entry, and the root node receives input from a spectrally positive L\'{e}vy process. Our aim is to characterize the distribution of the multivariate stationary workload under a specific scaling of the service rates. 
   Under mild conditions on the Laplace exponent of the driving L\'{e}vy process, we identify the limiting law of an appropriately scaled joint stationary workload in both light-traffic and heavy-traffic regimes. In particular, we establish conditions under which certain queueing workloads within the network asymptotically decouple, becoming independent in the limiting regime.

    \vspace{3mm}


\vspace{3mm}

\noindent
{\sc Affiliations.} NK is affiliated with the Korteweg-de Vries Institute for Mathematics, University of Amsterdam, Science Park 904, 1098 XH Amsterdam, The Netherlands.

\noindent
KD is affiliated with Mathematical Institute, Wroc\l{}aw University, pl.\ Grunwaldzki 2/4, 50-384 Wroc\l{}aw, Poland.

\noindent
MM is affiliated with the Mathematical Institute, Leiden University, P.O. Box 9512,
2300 RA Leiden,
The Netherlands. MM is also affiliated with Korteweg-de Vries Institute for Mathematics, University of Amsterdam, Amsterdam, The Netherlands; E{\sc urandom}, Eindhoven University of Technology, Eindhoven, The Netherlands; Amsterdam Business School, Faculty of Economics and Business, University of Amsterdam, Amsterdam, The Netherlands.

\vspace{3mm}

\noindent
{\sc Acknowledgments.}
KD is partially supported by the National Science Centre, Poland,  Grant No 2024/55/B/ST1/01062
(2025-2028). NK and MM are supported by the European Union’s Horizon 2020 research and innovation programme under the Marie Sk\l{}odowska-Curie grant agreement no.\ 945045, and by the NWO Gravitation project {\tiny NETWORKS} under grant agreement no.\ 024.002.003. \includegraphics[height=1em]{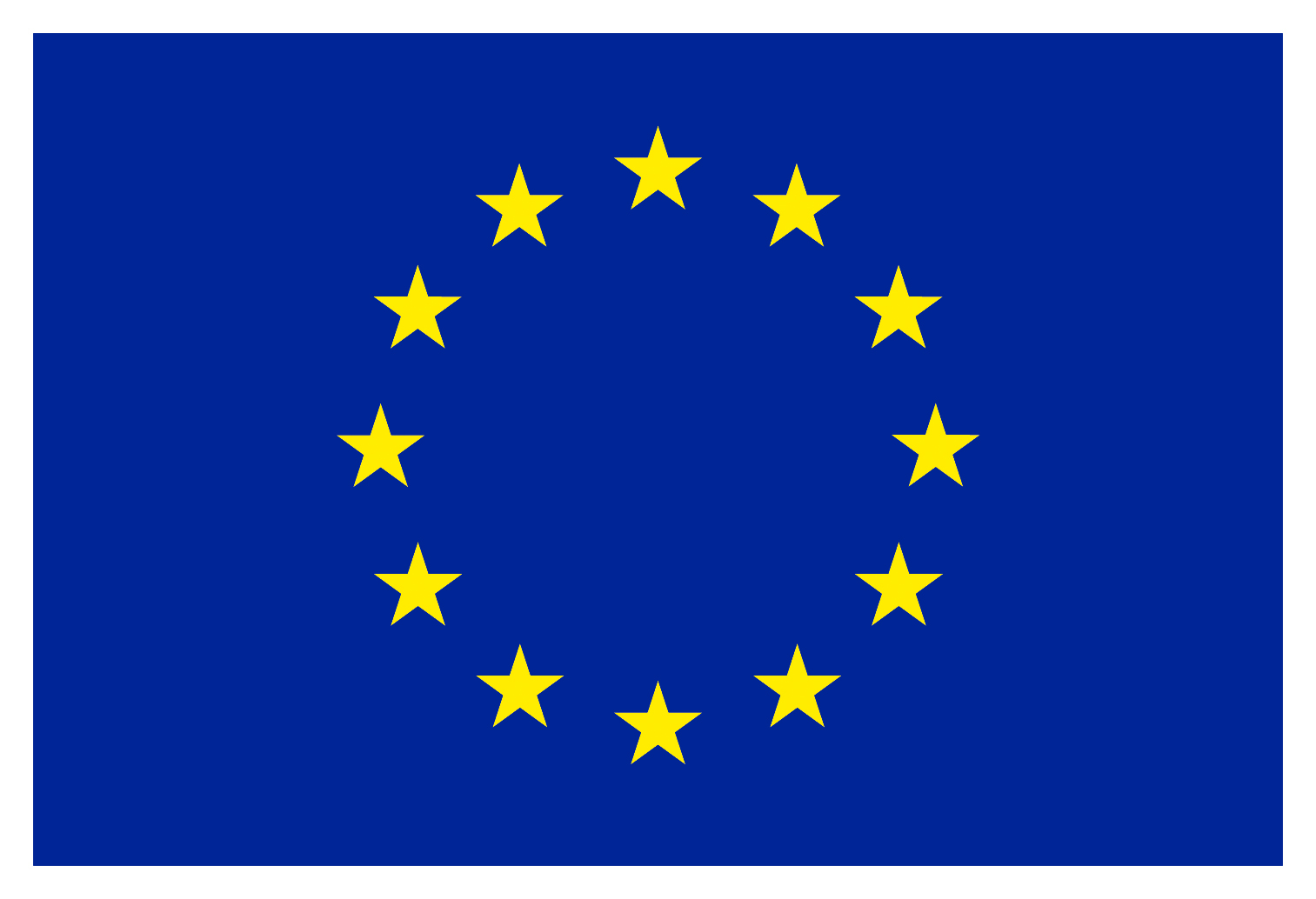}

\end{abstract}

\newpage

\section{Introduction}\label{sec:introduction}

We consider the classical model describing the evolution of a {\it storage system}, sometimes referred to as the {\it workload process} of a {\it queueing system}. The most basic model concerns a single buffered resource with unlimited storage capacity, fed by an external input process $J(\cdot)$, typically assumed to have stationary increments, and equipped with an output rate $r >0$. The storage level at time $t$ is a functional of the input process, and can be identified by solving the so-called {\it Skorokhod problem} \cite[\S IX.2]{Asmussen2003}. Indeed, it can be expressed in terms of the model primitives $J(\cdot)$ and $r$ via the representation
\begin{align}
    Q(t) = Q(0) + J(t) - rt +\max\biggl(0, \sup_{0<s<t}\bigl(-Q(0) - J(s)+rs\bigr)\biggr).\notag
\end{align}
Importantly, this formalism does not require that the process $J(\cdot)$ be increasing.
A vast branch of the literature aims to characterize the stationary workload $Q$, to be interpreted as the limit of $Q(t)$ as $t\to\infty$, for arbitrary initial workload $Q(0)$, under the stability condition ${\mathbb E}\{J(1)\} <r$. For a spectrally-positive L\'{e}vy input process $J(\cdot)$ (L\'evy processes without negative jumps, that is), the Laplace-Stieltjes transform (LST) pertaining to the stationary distribution is known from the (generalized version of the) Pollaczek-Khinchine theorem  \cite[\S III.1]{DeM15}.

While single queues have been very intensively studied, multi-queue networks have attracted significantly less attention. The simplest instance of such a network is a {\it tandem} queue, which consists $n\geqslant 2$ consecutively concatenated queues where the output of a queue serves as the input for the next one. We denote the respective output rates of these queues by $r_i$ for $i\in\{1,\ldots,n\}$. Without any loss of generality we can assume $r_1>\ldots>r_n$, as if for a pair of subsequent queues $i$ and $i+1$ this ordering would be violated, the $(i+1)$-st queue would remain empty and hence can be left out. Noting that the sum of the first $i$ workloads can be written as that of a single queue with output rate $r_i$, there is a way to describe the multivariate workload process $\vk Q(\cdot) = (Q_1(\cdot),\ldots, Q_n(\cdot))^{\top}$. More precisely, using that
\begin{align*}
\sum_{j=1}^i Q_i(t) = \left(\sum_{j=1}^{i}Q_i(0)\right) + J(t) - r_it +\max\left(0, \sup_{0<s<t}\left(-\left(\sum_{j=1}^{i}Q_i(0)\right) - J(s)+r_is\right)\right),
\end{align*}
we have, for any $t\geqslant 0$ and any $i\in\{2,\ldots,n\}$,
\begin{align*}
    Q_i(t)  = Q_i(0) - (r_i-r_{i-1})t &+ \max\left(0, \sup_{0<s<t}\left(-\left(\sum_{j=1}^{i}Q_i(0)\right) - J(s)+r_is\right)\right)\\
    & - \max\left(0, \sup_{0<s<t}\left(-\left(\sum_{j=1}^{i-1}Q_i(0)\right) - J(s)+r_{i-1}s\right)\right).
\end{align*}
Under the assumption  $\E{J(1)}<r_i$ the $i$-th queue is stable, so that the $n$-dimensional limiting workload ${\boldsymbol Q}$ is well defined.  For more background on multivariate storage processes, including the tandem systems discussed above, see e.g.\ \cite[Chs.\ XII-XIII]{DeM15} and \cite{Rob13}.
In the case that $J(\cdot)$ is a spectrally positive L\'{e}vy process, 
the joint stationary distribution pertaining to $\vk Q$ in the above class of tandem networks has been identified in terms of its LST;
see e.g.\ the results in  \cite{DDR}.

Besides tandem networks, one can consider more general classes of $n$-node {\it queueing networks}, where the subclass of {\it tree networks} still allows exact analysis. In this class of {\it feedforward queues}, work enters only at the root queue. At any node, the output of the queue can be \textit{split}, and either feeds into one or more subsequent queues or leaves the network. Importantly, this framework does not allow \textit{loops} or \textit{merges}; see again \cite{DDR}.
For the analysis of other extensions of the storage model with multidimensional L\'evy input  we refer to \cite{Kel1, Kel2, Kel3, Kel4, Kel5, Kon}.

We proceed by briefly introducing some notation; a more complete account of the model and the imposed assumptions is given in Section \ref{sec:model}. We denote by $p_{ij}\in[0,1]$ the fraction of the $i$-th node which serves as the input of $j$-th node.  
As we do not allow loops and merges in the network, the nodes can be ordered such that the matrix $P= \{p_{ij}\}_{i,j=1}^{n}$ 
is strictly upper-triangular and has at most one positive element per column. The model of a tandem queue is a particular example in which the entries of  $P$ are such that $p_{i,i+1}=1$  for $i\in\{1,\ldots,n-1\}$ and $p_{ij}=0$ otherwise. 
Hence, the queueing network that we consider can be described in terms of (i)~the {\it input process} $J(\cdot)$ (describing the work entering at queue $1$), (ii)~the {\it output rates} $r_1,\ldots,r_n$ of nodes $1,\ldots,n$, respectively, and (iii)~the {\it routing matrix} $P = \{p_{ij}\}_{i,j=1}^{n}$. Note that this queuing network model $(J(\cdot),P,{\vk r})$ covers tandem networks with `leaking', i.e., instances in which a part of the network's input leaves the system before queue $n$ has been reached (i.e., $p_{i,i+1}\in(0,1)$  for $i\in\{1,\ldots,n-1\}$ while still all other $p_{ij}=0$).

We center the input process $J(\cdot)$, which concretely entails that $\E{J(1)}=0$; observe that this can be done without losing any generality. 
In our setup the driving process $J(\cdot)$ is assumed to be a spectrally-positive L\'evy process. As was demonstrated in \cite{DDR}, the distribution of the stationary workload vector $\vk Q$ can be explicitly expressed in terms of $(J(\cdot),P,{\vk r})$ via its (joint) LST. We discuss this expression, along with the required assumptions, in Section \ref{sec:model}.

\medskip

\noindent {\it Scope \& contribution ---} 
The main objective of this paper is to analyze the tree network defined above in specific light-traffic and heavy-traffic regimes. As extensively motivated in, e.g., \cite[Ch.~XIV]{DeM15}, Lévy processes are particularly well suited to model aggregates of stochastic processes under appropriate time scalings, with applications in areas such as finance and communication networks. This perspective is further supported by the comprehensive treatment in \cite{Whitt2002},
and stochastic process limits arguments, as in \cite{MSR, TWS}. Additional background on the structure and properties of Lévy processes can be found in the classical monographs \cite{Sato1999,bertoin1996levy,Kypri}, while general heavy-traffic and functional limit principles for queueing networks are discussed in \cite{Harrison1985,ChenYao2001,Whitt2002}. The literature on queues in heavy traffic dates back to the seminal work of \cite{Kin}, while classical references on light-traffic behavior include \cite{DaR,Sig}.
In the setting of L\'evy inputs, the heavy traffic limit of a storage process was analyzed, e.g., in 
\cite{BKZ}, while the tandem system was investigated in \cite{Koo}.

More precisely, our goal is to describe the joint stationary distribution of all the nodes present in the network in the two following asymptotic regimes:
\begin{itemize}
\item[$\circ$] {\it light traffic}: for all nodes $i\in\{1,\ldots,n\}$ the output rates $r_i$ are such that $r_i=r_i(u)\to \infty$,
\item [$\circ$] {\it heavy traffic}: for all nodes $i\in\{1,\ldots,n\}$ the output rates $r_i$ are such that $r_i=r_i(u)\to 0$,
\end{itemize}
letting $u\to\infty$ (where $u$ denotes a scaling parameter). 
This paper contains two main results. In \netheo{main_light} we provide, under  mild assumptions on the network structure and the Laplace exponent of the input L\'{e}vy process, the limiting LST of the appropriately rescaled stationary workload in light traffic scenario. \netheo{main_heavy} is the counterpart of \netheo{main_light} for the heavy traffic scenario.
In particular, we provide the explicit conditions under which the stationary workloads for different queues `decouple', i.e., become independent in the limit as $u\to\infty$. Informally, these conditions entail that the workload processes at the different queues fluctuate at different timescales (`multiscale light- and heavy traffic'): a queue that fluctuates slowly sees a fast fluctuating queue essentially in equilibrium.

The asymptotic independence among stationary workloads under multiscale light- and heavy-traffic regimes has previously been established in various specific settings. In this context we refer to the seminal paper \cite{dai2023asymptotic} on the multiscale heavy-traffic limit theorem for Jackson networks in stationarity; see also the recent contributions \cite{dai2024asymptotic,guang2025steady}. In \cite{DKM}, the focus is on path-level convergence, with the main result providing functional limit theorems for tree networks with Gaussian inputs in multiscale light- and heavy-traffic. The approach taken in the present paper is intrinsically different from those in \cite{dai2023asymptotic, DKM} though: as we have access to the joint LST of the stationary workload vector, we can prove distributional convergence relying on L\'evy's convergence theorem. It is noted that one could alternatively first establish functional convergence and then analyze the corresponding stationary behavior. However, this approach would require addressing the delicate issue of interchanging the limit as time $t \to \infty$ with the limit of the scaling parameter $u \to \infty$.

\medskip

\noindent{\it Outline ---} This paper is organized as follows. {Subsection \ref{sub.21}} formally introduces the model and the assumptions imposed. In Subsection \ref{sec:notations} we provide all the notations which are necessary to present the main results. In Section \ref{sec:main} we present the two main results of this paper. Section \ref{sec:examples} contains a series of examples illustrating the application of the main results. This section has been split into two subsections: Subsection \ref{subsec:processes} presents examples of  L\'{e}vy processes which satisfy the imposed assumptions, while Subsection~\ref{subsec:networks} contains examples of particular networks. Section \ref{sec:proof} is dedicated to the proof of the main theorems.  Appendix~\ref{sec:appendix} contains (relatively technical) proofs of a series of auxiliary results.

\section{Model, assumptions and notation}\label{sec:model} 
This section details the model studied, the assumptions  imposed, and additional notations that are used throughout the paper. 

\subsection{Model and assumptions}\label{sub.21}
As discussed in the introduction, to describe the queueing network, we should specify the triplet $(J(\cdot),P,\vk r)$. We start by detailing the input process $J(\cdot)$, which serves as the unique input in the system, {feeding into node 1:} 
\begin{itemize}
    \item[\textbf{P}:]{$J(\cdot)$ is a centered spectrally positive L\'{e}vy process.}
\end{itemize}
Next, we specify the network {structure}, which involves the routing matrix $P= \{p_{ij}\}_{i,j=1}^{n}$ and the output rates $\vk r$. The output rates  $\vk r= \vk r(u)>0$ are parametrized by a scaling parameter $u$. To obtain a closed form of the stationary workload in terms of the driving process $J(\cdot)$, the associated Skorokhod problem needs to be solved; {see, e.g., \cite{DDR}}. To this end, we impose the following three assumptions:
\begin{itemize}
    \item[\textbf{N1}:]{The routing matrix $P$ is strictly upper triangular, contains only non-negative inputs, and each column, apart the first one, has exactly one positive entry;}
    \item[{\textbf{N2}$^\star$:}]{For any $i,j\in\{1,\ldots,n\}$, if $p_{ij}>0$ then $p_{ij}\,r_{i}(u)>r_j(u)$ for any $u>0$;}
    \item[{\textbf{N3}$^\star$:}]{For any $i\in\{1,\ldots,n\}$ and any $u>0$, {$\E{J(1)}\,(I-P^{\top})^{-1}_{i,1}<r_i(u)$}.}
\end{itemize}
Under {Assumptions \textbf{N1}, \textbf{N2}$^\star$ and \textbf{N3}$^\star$} the stationary workload can be expressed in terms of the input process $J(\cdot)$.
Indeed, as pointed out in e.g.\ \cite{DDR} and \cite[Ch.\ XIII]{DeM15}, under these conditions
the stationary workload $\vk Q_{u}$, for any given $u>0$ obeys the distributional equality
\begin{align}
    \vk Q_{u} \stackrel{\rm d}{=} (I - P^{\top})\overline{\vk X}_u,\label{W_def}
\end{align}
where $\overline{\vk X}_u = \bigl(\overline{X}_{u,1},\ldots,\overline{X}_{u,n}\bigr)^{\top}$ is
defined by
\begin{align*}
    \overline{\vk X}_u &:= \left(\sup_{{s\geqslant 0}} X_{u,1}(s),\,\ldots,\, \sup_{{s\geqslant 0}} X_{u,n}(s)\right)^\top,
\end{align*}
for the process $\vk X_u(\cdot) = (X_{1,u}(\cdot),\ldots,X_{n,u}(\cdot))$ defined for any $t\in\R$ as
\begin{align*}
    \vk X_u(t) &:= (I-P^{\top})^{-1}\vk J(t) - \vk r t = \hat{\vk p}\,J(t) - \vk r(u)\,t,
\end{align*}
where $\vk J(\cdot) := (J(\cdot),0,\ldots,0)\in\R^n$ and $\hat{\vk p}$ is the first column of the matrix $\left(I-P^\top\right)^{-1}$.

{In our approach, we intend to intensively work with LSTs, which means that we need to control the output rates for any pair of two nodes (i.e., not only for the consecutive nodes). As a consequence, Assumption \textbf{N2}$^\star$ does not suffice; instead, we need to assume that  the following more strict version of Assumption {\textbf{N2}$^{\star}$}} is in place:
\begin{itemize}
    \item[{\textbf{N2}}:]{For any $j\in\{1,\ldots, n-1\}$ and any $u>0$, we have $r_j(u)>0$ and $r_{j}(u)/\hat{p}_j> r_{j+1}(u)/\hat{p}_{j+1}$.}
\end{itemize}
{Indeed, the combination of Assumptions \textbf{N1} and \textbf{N2} implies that Assumption \textbf{N2}$^\star$ holds, which can be seen as follows.} For any $i,j\in\{1,\ldots n\}$ if $p_{ij}>0$, then according to Assumption \textbf{N1} we have that $j>i$, and from Assumption {\textbf{N2}} we have that
\[r_{i}(u)/\hat{p}_i > r_{i+1}(u)/\hat{p}_{i+1} >  \ldots > r_{j}(u)/\hat{p}_j.\] Hence, Assumption {\textbf{N2}$^\star$} follows from the fact that $\hat{p}_j = p_{ij}\hat{p}_i$ in case $p_{ij}>0$.

It can also be argued that Assumption \textbf{N3}$^\star$ immediately follows from Assumption \textbf{P}. Indeed, observe that $\E{J(1)}\,(I-P^{\top})^{-1}_{i,1}=0$ since we assumed that the process $J(\cdot)$ is centered. 
Considering the arguments presented in this and the preceding paragraphs, we will impose Assumptions $\textbf{P}$, $\textbf{N1}$ and $\textbf{N2}$ throughout this paper. Under these assumptions we know that $\textbf{N1}$, $\textbf{N2}^{\star}$ and $\textbf{N3}^\star$ are valid. 

Additionally, we impose the following restriction on the relative asymptotic behavior of the output rates:

\begin{itemize}
    \item[{\textbf{R}:}]{For any $i,j\in\{1,\ldots,n\}$ with $i>j$, we assume $\lim_{u\to\infty}r_i(u)/r_j(u)\in[0,1]$.}
\end{itemize}

In this paper we consider two limiting regimes: the light-traffic scenario, with $\vk r(u)\to \vk \infty$ as $u\to\infty$, and the heavy-traffic scenario, with $\vk r(u)\to \vk 0$ as $u\to\infty$ (where $\vk 0$ and $\vk\infty$ have the evident meaning). For each regime we impose specific assumptions on the Laplace exponent of the driving L\'evy process $J(\cdot)$:
\begin{itemize}
    \item[$\circ$]
For the light-traffic case we assume  
\begin{itemize}
        \item[\textbf{L}:]{$\log\E{e^{-sJ(1)}} = \mathfrak{C}s^{\alpha} + o(s^{\alpha})$ as $s\to\infty$ for some constants 
        {$\mathfrak{C}<0$} and $\alpha>1$.}
\end{itemize}
\item[$\circ$] For the heavy-traffic we assume  
\begin{itemize}
        \item[\textbf{H}:]{$\log\E{e^{-sJ(1)}} = \mathfrak{C}s^{\alpha} + o(s^{\alpha})$ as $s\downarrow 0$ for some constants {$\mathfrak{C}<0$} and $\alpha>1$.}
\end{itemize}
\end{itemize}

\subsection{Notation}\label{sec:notations}

Throughout the paper we consistently use the following notation. For any $n\in\N$, we write $\sprod{n}$ to denote the set $\{1,\ldots,n\}$, and for any $n,k\in\N$, $k\leqslant n$ we write $\sprod{n}^{k+}$ to denote the set $\{k,\ldots,n\}$ (so that {$\sprod{n}=\sprod{n}^{1+}$}). All vectors are written in bold and, vice versa, all bold symbols refer to vectors. 
If $\vk b\in\R^d$, we can also write $\vk b = (b_1,\ldots,b_d)^{\top}$ (i.e., for the entries of $\vk b$ we do not use the bold font); with a mild abuse of notation, if $\vk b\in\R^{\tt I}$ for some finite set ${\tt I} = \{i_1,\ldots,i_k\}\subset\N$, then $\vk b = (b_{i_1},\ldots,b_{i_n})^{\top}$. 
In this sense, $\R^{d}$ is a simplified notation for $\R^{\sprod{d}}$. For any finite set ${\tt I}$ we denote its cardinality by $\abs{{\tt I}}$.
{For any $\vk{a}, \vk{b} \in \R^k$,  we write $\vk a \circ \vk b$ to denote their component-wise multiplication, i.e.
$\vk a \circ \vk b:=(a_1b_1,\ldot,a_kb_k)^\top$.}

For any two functions $f(x),\,g(x)$, $x\in \R$ we write $f(x)\sim g(x)$ as $x\to\infty$ if $\lim_{x\to\infty}f(x)/g(x) = 1$. Likewise, we write $f(x) = o(g(x))$ as $x\to\infty$ if $\lim_{x\to\infty}f(x)/g(x) = 0$, and we write $f(x)\asymp g(x)$ as $x\to\infty$ if $\lim_{x\to\infty}f(x)/g(x)\in(0,\infty)$. 

{In addition to the general notation given above, we introduce further notation pertaining specifically to the network structure and the output rates. We start with the notation associated with the network structure.} For any node $j\in\sprod{n}^{2+}$ we denote its unique ancestor in the tree network by $j^{\prime}\in\sprod{n}$, i.e., the $j'$ such that $p_{j^{\prime}j}>0$. For any $j\in\sprod{n}$, the set ${\tt S}_j \subset\{j,\ldots,n\}$ is the set of vertices with index at least $j$, whose ancestor's index (if it exists) is smaller than $j$, i.e.,
\begin{align*}
    {\tt S}_j := \{j\}\cup\{i\in\sprod{n}^{(j+1)+}\colon i^{\prime}<j\}.
\end{align*}
Denote additionally for any $j\in\sprod{n}$ by ${\tt D}_j$ the set of all offsprings of the vertex $j$, i.e. 
\begin{align*}
    {\tt D}_j := \{i\in\sprod{n}\colon p_{ji}>0\}.
\end{align*}

We state below a number of relations between the sets ${\tt S}_j$ and ${\tt D}_{j}$ which we are going to repeatedly use in our proofs.

\begin{lem}\label{S_D_relations}
For any $j,k\in\sprod{n}$, $j\not = k$
\begin{align}
    &{\tt D}_j\cap {\tt D}_{k} = \varnothing,\label{D_intersection}\\
    &{\tt S}_{j}\cap\bigcup_{l=j}^{n}{\tt D}_{l} = \varnothing,\label{S_D_intersection}\\
    &{\tt S}_{j}\cup\bigcup_{l=j}^{n}{\tt D}_{l} = \sprod{n}^{j+},\label{S_D_union}\\
    &{\tt S}_{j+1} = {\tt S}_{j}\cup {\tt D}_{j}\setminus \{j\}.\label{S_rec}\\
    & k\in{\tt S}_{j} \iff j\in\{k^{\prime}+1, \ldots ,k\}.\label{k_in_S}
\end{align}
\end{lem}
The proof of \nelem{S_D_relations} can be found in Appendix \ref{A1}. In particular, \eqref{S_D_union} combined with \eqref{D_intersection} and \eqref{S_D_intersection} gives us that
\begin{align}
    {\tt S}_{j}\sqcup\bigsqcup_{l=j}^{n}{\tt D}_{l} = \sprod{n}^{j+},\label{S_D_strong_union}
\end{align}
where $A\sqcup B = C$ means that $A\cap B = \varnothing$ and $A\cup B = C$.
\medskip

{Next, we proceed to additional notation related to the output rates.} Using Assumption {\textbf{R}} one can divide the set $\sprod{n}$ into disjoint sets ${\tt I}_1,\ldots,{\tt I}_m$, for some $m\in\N$, and define a function $\mathcal{k}(\cdot): \sprod{n}\to\sprod{m}$ and, for any $k\in\sprod{m}$, a function $\mathcal{r}_k(\cdot)$ such that for any $i\in{\tt I}_k$
\begin{align}
    \lim_{u\to\infty}\frac{r_i(u)}{\mathcal{r}_{\mathcal{k}(i)}(u)} = \mathfrak{r}_i\in(0,1],\label{r_same}
\end{align}
and, for any $k,l\in\sprod{m}$ with $k>l$, 
\begin{align}
     \lim_{u\to\infty}\frac{\mathcal{r}_k(u)}{\mathcal{r}_l(u)} = 0.\label{r_diff}
\end{align}
Relation \eqref{r_diff} together with Assumption $\textbf{R1}$ tells us that the set ${\tt I}_k$ for any $k\in\sprod{m}$ is an `interval', i.e., if $a,b,c\in\sprod{n}$, $a\leqslant b\leqslant c$ and $a,c\in{\tt I}_k$, then $b\in{\tt I}_k$. For any {$k\in\sprod{m}$} we define the {minimal} element of the set ${\tt I}_k$ {by $q_k$, i.e.,}
\begin{align*}
    q_k := \min\{i\in{\tt I}_k\}.
\end{align*}
A crucial consequence of \eqref{r_same} and \eqref{r_diff} is that, for any $i,j\in\sprod{n}$ with $i\leqslant j$,
\begin{align}
    \lim_{u\to\infty}\frac{r_j(u)}{\mathcal{r}_{\mathcal{k}(i)}(u)} = 
        \mathfrak{r}_j\mathbb{I}_{\{j\in{\tt I}_{\mathcal{k}(i)}\}}\label{r_comp}
\end{align}

\medskip

Denote ${\tt S}^{\star}_{j} := {\tt S}_j\cap {\tt I}_{\mathcal{k}(j)}$ and ${\tt D}^{\star}_{j} := {\tt D}_j \cap {\tt I}_{\mathcal{k}(j)}$.

\begin{figure}

\tikzset{int/.style  = {draw, circle, node font=\small, minimum size=1cm}}
\tikzset{sq/.style  = {draw,rectangle, node font=\small, minimum height=1cm}}

\begin{tikzpicture}[auto,>=latex, scale=0.9]

\tikzstyle{sq}=[draw, rectangle, minimum width=1.3cm,
                minimum height=0.8cm, align=center, inner sep=2pt]
\tikzstyle{int}=[draw, circle, minimum size=0.8cm, align=center, inner sep=1pt]

\node[sq] (Q1) {$Q_1$};
\path (Q1.east) coordinate (aux1);
\node[int, anchor=west] (A) at (aux1) {\small $10u^2$};

\node (baux) [above right=1.0cm and 1.5cm of A] {};
\node[sq] (Q2) [right=0pt of baux] {$Q_2$};
\path[->] (A) edge node[above] {\small $\tfrac12$} (Q2);
\path (Q2.east) coordinate (aux2);
\node[int, anchor=west] (B) at (aux2) {\small $4u^2$};

\node (caux) [above right=1.0cm and 1.5cm of B] {};
\node[sq] (Q3) [right=0pt of caux] {$Q_3$};
\path[->] (B) edge node[above] {\small $\tfrac13$} (Q3);
\path (Q3.east) coordinate (aux3);
\node[int, anchor=west] (C) at (aux3) {\small $u^2$};

\path[->] (C) edge ++(1.5cm,0);

\node (daux) [right=1.5cm of B] {};
\node[sq] (Q4) [right=0pt of daux] {$Q_4$};
\path[->] (B) edge node[above] {\small $\tfrac13$} (Q4);
\path (Q4.east) coordinate (aux4);
\node[int, anchor=west] (D) at (aux4) {\small $2u$};

\path[->] (D) edge ++(1.5cm,0);

\node (eaux) [below right=1.0cm and 1.5cm of A] {};
\node[sq] (Q5) [right=0pt of eaux] {$Q_5$};
\path[->] (A) edge node[above] {\small $\tfrac12$} (Q5);
\path (Q5.east) coordinate (aux5);
\node[int, anchor=west] (E) at (aux5) {\small $4u$};

\path[->] (E) edge ++(1.5cm,0);

\node (faux) [below right=1.0cm and 1.5cm of B] {};
\node[sq] (Q6) [right=0pt of faux] {$Q_6$};
\path[->] (B) edge node[above] {\small $\tfrac13$} (Q6);
\path (Q6.east) coordinate (aux6);
\node[int, anchor=west] (F) at (aux6) {\small $u$};

\path[->] (F) edge ++(1.5cm,0);

\end{tikzpicture}

\caption{\label{F1}Example tree network with six nodes. The values of the $r_i(u)$ are given in the circles, the values of the $p_{ij}$ on top of the corresponding arrows. }
\end{figure}
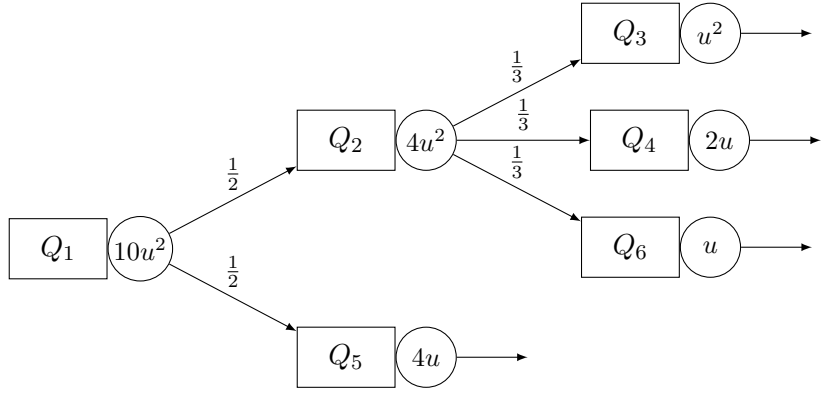

\medskip

\begin{example}{\em To illustrate the notation presented above we consider the following example network; see Figure \ref{F1}.
The transition matrix and output rate vector are given by, respectively,
\begin{align*}
P=\begin{pmatrix}
    0 & \tfrac 1 2 & 0 & 0 & \tfrac 1 2 & 0\\
    0 & 0 & \tfrac 1 3 & \tfrac 1 3 & 0 & \tfrac 1 3\\
    0 & 0 & 0 & 0 & 0 & 0\\
    0 & 0 & 0 & 0 & 0 & 0\\
    0 & 0 & 0 & 0 & 0 & 0\\
    0 & 0 & 0 & 0 & 0 & 0
\end{pmatrix},\quad
{\vk r}(u)=
\begin{pmatrix}
    10u^2 \\
    4u^2\\
    u^2\\
    2u\\
    4u\\
    u 
\end{pmatrix}.
\end{align*}
It is now readily verified that $\hat{\vk p} = (1,\,1/2,\,1/6,\,1/6,\,1/2,\,1/6)^{\top}$. This means that Assumption \textbf{N2} is satisfied, in that for $u>1$ we have
\begin{align*}
    r_1(u)\,>\,2r_2(u)\,>\,6r_3(u)\,>\,6r_4(u)\,>\,2r_5(u)\,>\,6r_6(u).
\end{align*}
The output rates have only to equivalence classes, namely
\begin{align*}
    {\tt I}_1 = \{1,2,3\},\qquad {\tt I}_2 = \{4,5,6\},
\end{align*}
while it can be checked that the sets ${\tt S}_j$, ${\tt D}_j$, ${\tt S}^{\star}_j$ and ${\tt D}^{\star}_j$ for $j\in\sprod{6}$ are 
\begin{align*}
    {\tt S}_1 &= \{1\},\quad {\tt S}^{\star}_1 = \{1\},\quad {\tt D}_1 = \{2,3\},\quad {\tt D}^{\star}_1 = \{2\},\\
    {\tt S}_2 &= \{2,5\},\quad {\tt S}^{\star}_2 = \{2\},\quad {\tt D}_2 = \{3,4,6\},\quad {\tt D}^{\star}_2 = \{3\},\\
    {\tt S}_3 &= \{3,4,5,6\},\quad {\tt S}^{\star}_3 = \{3\},\quad {\tt D}_3 = \{\varnothing\},\quad {\tt D}^{\star}_3 = \{\varnothing\},\\
    {\tt S}_4 &= \{4,5,6\},\quad {\tt S}^{\star}_4 = \{4,5,6\},\quad {\tt D}_4 = \{\varnothing\},\quad {\tt D}^{\star}_4 = \{\varnothing\},\\
    {\tt S}_5 &= \{5,6\},\quad {\tt S}^{\star}_5 = \{5,6\},\quad {\tt D}_5 = \{\varnothing\},\quad {\tt D}^{\star}_5 = \{\varnothing\},\\
    {\tt S}_6 &= \{6\},\quad {\tt S}^{\star}_6 = \{6\},\quad {\tt D}_6 = \{\varnothing\},\quad {\tt D}^{\star}_6 = \{\varnothing\}.
\end{align*}
The above example illustrates how to identify the sets that we introduced above. \hfill $\Diamond$}
\end{example}

To state the {main results of this contribution}, we need some additional notation.
Let  $\Psi_{\alpha,C,j}^{-1}(\cdot)$ be (for $C>0$, $\alpha>1$ and $j\in\sprod{n}$) the inverse of the function $\Psi_{\alpha,C,j}(\cdot)$ defined by
\begin{align}
     \Psi_{\alpha,C,j}(s) &:= \mathfrak{r}_js + C\hat{p}_j^{\alpha}s^{\alpha}\label{Psi_def}.
     \end{align}
Then    define (for $k\in\sprod{m}$,  $j\in\sprod{n}$, $C>0$, $\alpha>1$ and any $\vk\omega\in\vk\R^{{\tt I}_{\mathcal{k}(j)}}$) the following constants:
\begin{align}
    \mathcal{A}_{\alpha,k}(C,\vk\omega) &:= \sum\limits_{l\in{\tt I}_{k}}\frac{\mathfrak{r}_{l}}{\hat{p}_{l}}\sum\limits_{m\in{\tt D}^{\star}_l}\hat{p}_m\omega_m - \sum\limits_{l\in{\tt I}_{k}}\omega_l\mathfrak{r}_l - C \left(\sum\limits_{l\in{\tt S}^{\star}_{q_k}}\hat{p}_{l}\omega_l\right)^{\alpha},\label{A_cal}\\
    \mathcal{C}_{\alpha,j}(C,\vk\omega) &:=  \Psi_{\alpha,C,j}^{-1}\bigl(\mathcal{F}_{j}(\vk\omega)\bigr)-\sum_{l\in{\tt S}^{\star}_j}\frac{\hat{p}_{l}\omega_l}{\hat{p}_j},\label{C_cal}\\
    \mathcal{D}_{\alpha,j}(C,\vk\omega) &:=  \Psi_{\alpha,C,j}^{-1}\bigl(\mathcal{F}_{j}(\vk\omega)\bigr) -\sum_{l\in{\tt S}^{\star}_{j+1}}\frac{\hat{p}_{l}\omega_l}{\hat{p}_j},\label{D_cal}
\end{align}
where 
\begin{align}
 \mathcal{F}_{j}(\vk\omega) &:= \sum_{\substack{l=j+1 \\ l\in{\tt I}_{\mathcal{k}(j+1)}}}^{n}\left(\frac{\mathfrak{r}_{l-1}}{\hat{p}_{l-1}} - \frac{\mathfrak{r}_l}{\hat{p}_{l}}\right)\sum_{i\in{\tt S}^{\star}_{l}}\hat{p}_{i}\omega_i.\label{F_cal}
\end{align}

\COM{\begin{remark} The form of the inverse functions $\Psi^{-1}_{\alpha,C,j}(\cdot)$ can be unified over all $j\in\sprod{n}$. The equation (for $s\in\R_+^{n}$)
    \begin{align*}
        s = \mathfrak{r}_j\Psi^{-1}_{\alpha,C,j}(s) + C\hat{p}_j^{\alpha}\left(\Psi^{-1}_{\alpha,C,j}(s)\right)^{\alpha}
    \end{align*}
    is equivalent to
    \begin{align*}
        a_j s &= a_j\mathfrak{r}_j\Psi^{-1}_{\alpha,C,j}(s) + Ca_j\hat{p}_j^{\alpha}\left(\Psi^{-1}_{\alpha,C,j}(s)\right)^{\alpha}\\
        &= a_j\mathfrak{r}_j\Psi^{-1}_{\alpha,C,j}(s) +C\left(a_j\mathfrak{r}_j\Psi^{-1}_{\alpha,C,j}(s)\right)^{\alpha},
    \end{align*}
    where (denote $\beta := 1/(\alpha-1)$)
    \begin{align*}
        a_j = \hat{p}^{\alpha\beta}_j\mathfrak{r}^{-\alpha\beta}_{j}.
    \end{align*}
    Hence, if we denote $\widetilde{\Psi}_{\alpha,C}(s) = s+Cs^{\alpha}$ then, for any $s\in\R_+$
    \begin{align*}
        a_j\mathfrak{r}_j\Psi^{-1}_{\alpha,C,j}(s) = \widetilde{\Psi}_{\alpha,C}^{-1}(a_j s),
    \end{align*}
    implying that
    \begin{align*}
        \Psi^{-1}_{\alpha,C,j}(s) = \hat{p}^{-\alpha\beta}_j\mathfrak{r}^{\beta}_{j}\widetilde{\Psi}_{\alpha,C}^{-1}\left(\hat{p}^{\alpha\beta}_j\mathfrak{r}^{-\alpha\beta}_{j}s\right).
    \end{align*}
\end{remark}}

\section{Main result} \label{sec:main}
In this section we present the main results of this paper.  Our goal is to describe, as $u\to\infty$, the joint distribution of the random variables $Q^{r,u}_j:=r_j^{\beta}(u)\,Q_{u,j}$, for $j\in\sprod{n}$, denoting $\beta := 1/(\alpha-1)$). 
We separate our results into two theorems: \netheo{main_light} covers the light-traffic case, while \netheo{main_heavy} the heavy-traffic case. Recall that we have split $\sprod{n}$ into the $m$ disjoint sets ${\tt I}_1,\ldots,{\tt I}_m$; see Subsection \ref{sec:notations}.

\begin{theo}\label{main_light} Let the network $(J(\cdot),\,P,\,\vk r(\cdot))$ be in the light-traffic scenario. Grant Assumptions \textbf{\rm \bf P}, \textbf{\rm \bf N1}, \textbf{\rm \bf N2}, \textbf{\rm \bf R} and \textbf{\rm \bf L}. Then, for any $\vk \omega\in \R_+^{n}$,
    \begin{align*}
        \lim_{u\to\infty}\E{e^{-\sprod{\vk\omega, \vk{Q}^{r,u}}}} = \prod_{k=1}^{m}F_k\left(\vk{\mathfrak{r}}^{1/(\alpha-1)}_{{\tt I}_k}\circ\vk\omega_{{\tt I}_k}\right),
    \end{align*}
    where, for any $\vk \omega\in \R_{+}^{{\tt I}_k}$, 
    \begin{align*}
        F_k(\vk\omega) := \frac{\omega_{q_k+\abs{{\tt I}_{k}}-1}\mathfrak{r}_{q_k+\abs{{\tt I}_{k}}-1}}{\abs{\mathcal{A}_{\alpha,k}(\mathfrak{C},\vk\omega)}}\prod_{j=0}^{\abs{{\tt I}_k}-2}\frac{\abs{\mathcal{C}_{\alpha,q_k+j}(\mathfrak{C},\vk\omega)}}{\abs{\mathcal{D}_{\alpha,q_k+j}(\mathfrak{C},\vk\omega)}}.
\end{align*}
\end{theo}

The proof of \netheo{main_heavy} is presented in Section \ref{sec:proof}.

\begin{theo}\label{main_heavy} Let the network $(J(\cdot),\,P,\,\vk r(\cdot))$ be in the heavy-traffic scenario. Grant Assumptions \textbf{\rm \bf P}, \textbf{\rm \bf N1}, {\textbf{\rm \bf N2}, \textbf{\rm \bf R}} and \textbf{\rm \bf H}. Then, for any $\vk \omega\in \R_+^{n}$,
    \begin{align*}
        \lim_{u\to\infty}\E{e^{-\sprod{\vk\omega, \vk{Q}^{r,u}}}} = \prod_{k=1}^{m}F_k\left(\vk{\mathfrak{r}}^{1/(\alpha-1)}_{{\tt I}_k}\circ\vk\omega_{{\tt I}_k}\right),
    \end{align*}
    where, for any $\vk \omega\in \R_{+}^{{\tt I}_k}$, 
    \begin{align*}
        F_k(\vk\omega) := \frac{\omega_{q_k+\abs{{\tt I}_{k}}-1}\mathfrak{r}_{q_k+\abs{{\tt I}_{k}}-1}}{\abs{\mathcal{A}_{\alpha,k}(\mathfrak{C},\vk\omega)}}\prod_{j=0}^{\abs{{\tt I}_k}-2}\frac{\abs{\mathcal{C}_{\alpha,q_k+j}(\mathfrak{C},\vk\omega)}}{\abs{\mathcal{D}_{\alpha,q_k+j}(\mathfrak{C},\vk\omega)}}.
\end{align*}
\end{theo}

The proof of \netheo{main_heavy} is analogous to the proof of \netheo{main_light}.

\begin{remark} {\em 
    Importantly, \netheo{main_heavy} requires an assumption on the asymptotics of the Laplace exponent of $J(\cdot)$ at zero (namely $\mathbf{L}$), whereas \netheo{main_light} requires an assumption on the asymptotics of the Laplace exponent of $J(\cdot)$ at  infinity (namely $\mathbf{H}$). Section \ref{sec:examples} indicates that Assumption \textbf{H} is satisfied by a broad class of L\'{e}vy processes, while a considerably more narrow class satisfies Assumption \textbf{L}. }\hfill $\Diamond$
\end{remark}

\begin{remark}\label{rem:zero_D} {\em 
    Our assumptions cannot guarantee the sign of the constants $\mathcal{A}_{\alpha,k}(\mathfrak{C},\vk\omega)$ and $\mathcal{D}_{\alpha,q_k+j}(\mathfrak{C},\vk\omega)$, or even that these constants are non-zero. However, if $\mathcal{A}_{\alpha,k}(\mathfrak{C},\vk\omega^{\star})=0$ for some $\vk\omega^{\star}\in\R_+^n$, then this implies that 
    \begin{align}
        \lim_{\omega\to\omega^{\star}} \mathcal{C}_{\alpha,q_k}(\mathfrak{C},\vk\omega)/\mathcal{A}_{\alpha,k}(\mathfrak{C},\vk\omega)\in(0,\infty);\label{rem_1}
    \end{align}
    if $\mathcal{D}_{\alpha,q_k+j}(\mathfrak{C},\vk\omega^{\star}) = 0$ for some $\vk\omega\in\R_+^{\abs{{\tt I}_k}}$, then this  implies that $j<\abs{{\tt I}_k}-2$, and 
    \begin{align}
        \lim_{\omega\to\omega^{\star}}  \mathcal{C}_{\alpha,q_k+j+1}(\mathfrak{C},\vk\omega)/\mathcal{D}_{\alpha,q_k+j}(\mathfrak{C},\vk\omega)\in(0,\infty).\label{rem_2}
    \end{align}
    Hence, the functions $F_k(\vk\omega)$ can still be well-defined for all $\vk\omega\in\R_+^{{\tt I}_k}$. We relegate the proof of these properties to Appendix \ref{A3}. }\hfill $\Diamond$
\end{remark}

\section{Examples}\label{sec:examples}
We can separate the assumptions imposed in Theorems \ref{main_light} and \ref{main_heavy} into two categories: assumptions on the input process $J(\cdot)$ (i.e., \textbf{P}, \textbf{L}, \textbf{H}) and assumptions on the network structure (i.e., \textbf{N1}, {\textbf{N2} and \textbf{R}}). In this section we first present two examples of input process covered by our results, and after that we discuss the application of Theorems \ref{main_light} and \ref{main_heavy} for two specific networks (with the required assumptions on the process $J(\cdot)$ being in place). 

\subsection{Examples of the input process}\label{subsec:processes} In this subsection we discuss three classes of processes: compound Poisson processes, Gamma processes, and sums of $\alpha$-stable processes.

\medskip 

{\it Centered spectrally-positive compound Poisson processes --- } Let $J^\star(t)$ for $t\in[0,\infty)$ be a compound Poisson process with intensity $\lambda>0$ and with independent almost-sure non-negative jumps that are distributed as the random variable $B$ with LST $b(s) := \E{e^{-sB}}$. We denote by $J(\cdot)$ the centered version of the process $J^\star(\cdot)$, i.e., for any $t\in[0,\infty)$,
\begin{align*}
    J(t) := J^\star(t) - \E{J^\star(t)} = J^\star(t) - \lambda t \,\E{B}.
\end{align*}
The process $J(\cdot)$ has only positive jumps, so that it satisfies \textbf{P}. The Laplace expoenent pertaining to the process $J(\cdot)$ is, for any $t\in[0,\infty)$,
\begin{align*}
    \log \E{e^{-sJ(1)}} = \log \E{e^{-sJ^\star(1)}}+{\lambda s \, \E{B}} = {\lambda \,\big(b(s) - 1 + s\,\E{B}\big)}.
    \end{align*}
In case $\E{B^2}<\infty$, we can write \[\log \E{e^{sJ(1)}} = s^2\lambda\E{B^2}/2 + o(s^2)\]
as $s\downarrow 0$, 
implying that process $J(\cdot)$ satisfies Assumption \textbf{H} with $\alpha = 2$ and $\mathfrak{C} = \lambda\E{B^2}/2$.
However, in the domain $s\to\infty$ we have that $b(s)\to 0$, implying that $\log \E{e^{-sJ(1)}} = \E{B}s + o(s)$ as $s\to\infty$, so that this class of processes cannot satisfy Assumption \textbf{L}.

\medskip

{\it Centered Gamma processes --- } For positive constants $\gamma,\beta>0$, consider the L\'{e}vy process $J(\cdot)$ characterized through its {\it L\'evy triple} $(d,0,\Pi)$ \cite[Section 1.1, Subsection 1.2.4]{Kypri}, where
\begin{align*}
    \log\E{e^{-sJ(1)}} =-s d + \beta \log\left(\frac{\gamma}{\gamma+{s}}\right)+s \int_0^1 x\,\Pi(\td x), \:\:\mbox{with}\:\: \Pi(\td x) := \mathbb{I}_{x>0}\frac{\beta}{x}e^{-\gamma x}\td x.
\end{align*}
As the support of $\Pi$ lies exclusively on $\R_+$, the process $J(\cdot)$ is spectrally positive. Noting that
\[\E{J(1)}=-\lim_{s\downarrow 0}\frac{{\rm d}}{{\rm d}s} \log\E{e^{-sJ(1)}}= d+\frac{\beta}{\gamma}-\int_0^1 x\,\Pi(\td x), \]
the process is centered if we pick 
\[d:= \int_0^1 x\Pi({\rm d}x)- \frac{\beta}{\gamma}.\]
With this choice, we conclude that $J(\cdot)$ satisfies Assumption \textbf{P}. Applying the Frullani integral, we obtain  
\begin{align*}
    \E{e^{-sJ(1)}} = \left(\frac{\gamma}{\gamma + s}\right)^{\beta}e^{{s\beta}/{\gamma}}; 
\end{align*}
see e.g.\ \cite[Section 2.2]{DeM15}.
Hence, as $s\downarrow 0$,
\begin{align*}
    \log \E{e^{-sJ(1)}} &= \beta\log\left(\frac{\gamma}{\gamma + s}\right) + \frac{s\beta}{\gamma} = \frac{s\beta}{\gamma} - \frac{s\beta}{\gamma + s} - \frac{s^2\beta}{2(\gamma+s)^2} +o(s^2)\\
    &= \frac{s^2\beta}{\gamma(\gamma+s)} - \frac{s^2\beta}{2(\gamma+s)^2}+o(s^2) = \frac{s^2\beta}{2\gamma^2} + o(s^2),
\end{align*}
implying that $J(\cdot)$ satisfies Assumption \textbf{H} with $\alpha = 2$ and $\mathfrak{C} = \beta/2\gamma^2$. However, as $s\to\infty$ we have that $\log \E{e^{-sJ(1)}} = s\beta/\gamma + o(s)$, so that the centered Gamma process $J(\cdot)$ does not satisfy Assumption \textbf{L}.

\medskip

{\it {Superposition} of $\alpha$-stable processes --- } Fix some integer constant $l\in\N$. Let $J_k(\cdot)$ for any $k\in\sprod{l}$ be a centered, spectrally-positive $\alpha$-stable L\'{e}vy process with parameter $\alpha_k\in(1,2]$. We assume that all the processes $J_k(\cdot)$ for $k\in\sprod{l}$ are mutually independent. According to \cite[Section 2.3]{DeM15}, the Laplace exponent corresponding to $J_k(1)$ has the following form: for a constant $C_k$,
\begin{align*}
    \log \E{e^{-sJ_k(1)}} = {C_k s^{\alpha_k}}.
\end{align*}
Then the process $\smash{J(\cdot) := \sum_{k=1}^{l}J_k(\cdot)}$ is still {a} centered spectrally-positive L\'{e}vy process, hence it satisfies Assumption \textbf{P}, and the Laplace exponent of $J(\cdot)$ is
\begin{align*}
    \log\E{e^{-sJ(1)}} = \sum_{k=1}^{l} C_k s^{\alpha_k}.
\end{align*}
Hence, the process $J(\cdot)$ satisfies Assumption \textbf{L} with
\begin{align*}
    &\alpha = \max_{k\in\sprod{l}}\alpha_k,\qquad \mathfrak{C} = \sum_{\substack{k=1}}^{l}C_k \mathbb{I}_{\{\alpha_k = \alpha\}},
\end{align*}
and it satisfies Assumption \textbf{H} with
\begin{align*}
    &\alpha = \min_{k\in\sprod{l}}\alpha_k,\qquad \mathfrak{C} = \sum_{\substack{k=1 }}^{l}C_k\mathbb{I}_{\{\alpha_k = \alpha\}}.
\end{align*}

Similarly, we can consider various mixtures of the three types of processes discussed above, such as a convolution of a spectrally-positive $\alpha$-stable L\'{e}vy process with a spectrally positive compound Poisson, or with a Gamma process.

\subsection{Examples of the network structure}\label{subsec:networks} In this section we illustrate our main results by considering specific network structures.

\medskip 

{\it Two-layers network --- }
We consider a model with two layers: one root node and $n-1$ nodes in the second layer, characterized via the routing matrix
\begin{align*}
    P=\begin{pmatrix} 0 & p_2 &\ldots & p_n \\ 0 & 0 &\ldots & 0 \\ \vdots & \vdots &\ddots & \vdots \\ 0 & 0 &\ldots & 0  \end{pmatrix}.
\end{align*}
The output rates $r_i(u)$ are such that for any $i\in\sprod{n}^{2+}$
\begin{align*}
    \lim_{u\to\infty}\frac{r_i(u)}{r_1(u)} = 0,\quad 
    \lim_{u\to\infty}\frac{r_i(u)}{r_2(u)} = \mathfrak{r}_i.
\end{align*}
for some constants $\mathfrak{r}_j > 0$. Assume additionally that for any $j\in\sprod{n-1}^{2+}$ and any $u>0$
\begin{align*}
    \frac{r_{j+1}(u)}{p_{j+1}} < \frac{r_j(u)}{p_j}.
\end{align*}
Then, if the input process $J(\cdot)$ satisfies Assumptions \textbf{P} and \textbf{L} in the case of light traffic, or Assumptions \textbf{P} and \textbf{H} in the case of heavy traffic,  our results entail that
\begin{align}\label{tree}
   \lim_{u\to\infty}\E{e^{-\sprod{\vk\omega, \vk{Q}^{r,u}}}} = F_1(\omega_1)F_2(\mathfrak{r}^{\beta}_2\omega_2,\ldots,\mathfrak{r}^{\beta}_{n}\omega_n)
\end{align}
where, with $\psi_j(s) = \mathfrak{r}_js + \mathfrak{C}p_j^{\alpha}s^{\alpha}$ for $\alpha\in(1,2]$,
\begin{align*}
    F_1(\omega_1) &= \frac{1}{1 + \mathfrak{C}\omega_1^{\alpha-1}},\quad
   F_2(\vk\omega) = \frac{\omega_{n}\mathfrak{r}_{n}}{\sum\limits_{l=2}^{n}\omega_l\mathfrak{r}_l + \mathfrak{C} \left(\sum\limits_{l=2}^{n}p_l\omega_l\right)^{\alpha}}\prod_{j=2}^{n-1}\frac{\abs{\psi_{j}^{-1}\left(\sum\limits_{l=j+1}^{n} \left(\frac{\mathfrak{r}_j}{p_j} - \frac{\mathfrak{r}_l}{p_l}\right)p_l\omega_l\right)-\sum\limits_{l = j}^{n}\frac{p_{l}\omega_l}{p_j}}}{\abs{\psi_{j}^{-1}\left(\sum\limits_{l=j+1}^{n} \left(\frac{\mathfrak{r}_j}{p_j} - \frac{\mathfrak{r}_l}{p_l}\right)p_l\omega_l\right)-\sum\limits_{l = j+1}^{n}\frac{p_{l}\omega_l}{p_j}}}.
\end{align*}
The proof underlying \eqref{tree} is given in Appendix \ref{A4}.
Observe that in the limiting regime considered, there is asymptotic independence between the stationary workload in the root node and the stationary workloads in the downstream nodes, but the individual downstream nodes remain dependent; cf.\ \cite{dai2023asymptotic,DKM}. It is seen that the stationary workload in the root queue has a limiting distribution of the {\it Mittag-Leffler} type  (of index $\alpha$). 

We proceed by detailing a few  particular cases. In the first place, if $\mathfrak{r}_2=\ldots = \mathfrak{r}_{n} = 1$, then $\psi_j(s) = s+\mathfrak{C}p_j^{\alpha}s^{\alpha}$. We find that 
\begin{align*}
    F_2(\omega_2,\ldots,\omega_n) = \frac{\omega_n}{\sum\limits_{l=2}^{n}\omega_l + \mathfrak{C} \left(\sum\limits_{l=2}^{n}p_{l}\omega_l\right)^{\alpha}}\prod_{j=2}^{n-1}\frac{\psi_j^{-1}\left( \sum\limits_{l=j+1}^{n}\frac{p_l-p_j}{p_j}\omega_l\right) -\sum\limits_{l=j}^{n}\frac{p_{l}\omega_l}{p_j}}{\psi_j^{-1}\left( \sum\limits_{l=j+1}^{n}\frac{p_l-p_j}{p_j}\omega_l\right) -\sum\limits_{l=j+1}^{n}\frac{p_{l}\omega_l}{p_j}}.
\end{align*}
Second, if $p_2 = \ldots = p_n = 1/(n-1)$, then $\psi_j(s) = \mathfrak{r}_j s + \frac{\mathfrak{C}}{(n-1)^{\alpha}}s^{\alpha}$. It now follows that
\begin{align*}
   F_2(\omega_2,\ldots,\omega_n) = \frac{\omega_n\mathfrak{r}_n}{ \sum\limits_{l=2}^{n}\omega_l\mathfrak{r}_{l} + \mathfrak{C} \left(\sum\limits_{l=2}^{n}\frac{\omega_l}{n-1}\right)^{\alpha}}\prod_{j=2}^{n-1}\frac{\psi_j^{-1}\left(\sum\limits_{l=j+1}^{n}(\mathfrak{r}_{j} - \mathfrak{r}_l)\omega_l\right) -\sum_{l=j}^{n}\omega_l\mathfrak{r}_l}{\psi_j^{-1}\left(\sum\limits_{l=j+1}^{n}(\mathfrak{r}_{j} - \mathfrak{r}_l)\omega_l\right) -\sum_{l=j+1}^{n}\omega_l\mathfrak{r}_l}.
\end{align*}
If both $\mathfrak{r}_2=\ldots = \mathfrak{r}_{n} = 1$ and $p_2 = \ldots = p_n = 1/(n-1)$, then
\begin{align*}
    F_2(\omega_2,\ldots,\omega_n) = \frac{\sum\limits_{l=2}^{n}\omega_l}{\sum\limits_{l=2}^{n}\omega_l - \mathfrak{C} \,\left(\frac{1}{n-1}\sum\limits_{l=2}^{n}\omega_l\right)^{\alpha}}= \frac{1}{1-{\displaystyle \frac{{\mathfrak C}}{(n-1)^\alpha}}\left(\sum\limits_{l=2}^{n}\omega_l\right)^{\alpha-1}},
\end{align*}
which entails that in the limiting regime considered the sum of the $n-1$ downstream queues, in stationarity, follows the {Mittag-Leffler distribution} (of index $\alpha$).

\medskip

{\it Tandem --- }
As a second example, we consider an $n$-nodes tandem system, i.e., a network characterized via the routing matrix 
\begin{align*}
    P=\begin{pmatrix} 0 & 1 & 0 &\ldots & 0 \\  0 & 0 & 1 &\ldots & 0  \\ \vdots & \vdots & \vdots & \ddots & \vdots \\  0 & 0 & 0 &\ldots & 1 \\ 0 & 0 & 0 &\ldots & 0  \end{pmatrix},
\end{align*}
with output rates $r_1(u)>\ldots > r_n(u)$. {We assume, for any $i,j\in\sprod{n}$ with $i<j$, that}
\begin{align*}
    {\lim_{u\to\infty}\frac{r_j(u)}{r_i(u)}\in[0,1].}
\end{align*}
If the input process $J(\cdot)$ satisfies Assumptions \textbf{P} and \textbf{L} in case of light traffic, or Assumptions \textbf{P} and \textbf{H} in case of heavy traffic, then
\begin{align}
   {\lim_{u\to\infty}\E{e^{-\sprod{\vk\omega, \vk{Q}^{r.u}}}} = \prod_{k=1}^{m}F_k(\omega_{q_k},\ldots,\omega_{q_{k+1}-1}),}\label{tandem}
\end{align}
{where $q_{m+1} = n+1$, and the integer-valued sequence $1=q_1<\ldots < q_m\leq n$ is chosen in such a way that for any $k\in\sprod{m}$ and any $j\in\sprod{q_{k+1}-1}^{q_k+}$ (denoting $r_{n+1}(u):=0$)}
\begin{align*}
    &{\lim_{u\to\infty}\frac{r_{q_{k+1}}(u)}{r_{q_{k}}(u)}=0,} \quad {\lim_{u\to\infty}\frac{r_{j}(u)}{r_{q_{k}}(u)}=\mathfrak{r}_j>0,}
\end{align*}
(noting that such a sequence exists and  unique), {with} $\psi_j(s) := \mathfrak{r}_j s+\mathfrak{C}s^{\alpha}$ and $F_k\equiv F_k(\omega_{q_k},\ldots,\omega_{q_{k+1}-1})$ given by
\begin{align*}
    {F_k:=\abs{\frac{\omega_{q_{k+1}-1}\mathfrak{r}^{\alpha\beta}_{q_{k+1}-1}}{\sum\limits_{l=q_k+1}^{q_{k+1}-1}\bigl(\mathfrak{r}_{l-1} - \mathfrak{r}_{l}\bigr)\omega_{l}\mathfrak{r}_{l}^{\beta} - \omega_{q_k} - \mathfrak{C} \omega_{q_k}^{\alpha}}\prod_{j=q_k}^{q_{k+1}-2}\frac{\psi_j^{-1}\left(\sum\limits_{l=j+1}^{q_{k+1}-1}\bigl(\mathfrak{r}_{l-1} - \mathfrak{r}_{l}\bigr)\omega_{l}\mathfrak{r}_{l}^{\beta}\right) -\omega_j\mathfrak{r}^{\beta}_j}{\psi_j^{-1}\left(\sum\limits_{l=j+1}^{q_{k+1}-1}\bigl(\mathfrak{r}_{l-1} - \mathfrak{r}_{l}\bigr)\omega_{l}\mathfrak{r}_{l}^{\beta}\right) -\omega_{j+1}\mathfrak{r}^{\beta}_{j+1}}}.}
\end{align*}
The proof underlying \eqref{tandem} is given in Appendix \ref{A5}. In particular, in the case of `complete decoupling' (i.e., if $\lim_{u\to\infty} r_j(u)/r_i(u) =0$ for any $i,j\in\sprod{n}$ with $i<j$), we have $m=n$, $q_k=k$ for any $k\in\sprod{n}$, and
\begin{align*}
    {\lim_{u\to\infty}\E{e^{-\sprod{\vk\omega, \vk{Q}^{r,u}}}} = \prod_{k=1}^{m}\frac{1}{1 + \mathfrak{C} \omega_{k}^{\alpha-1}}},
\end{align*}
{implying that in this limiting regime the stationary workloads of all the nodes are asymptotically independent and identically distributed, with a limiting distribution that is of the Mittag-Leffler type of index $\alpha$.}

\section{Proof of \netheo{main_light}}\label{sec:proof}

Throughout this proof we consider all asymptotic relations as $u\to\infty$, unless specificly mentioned otherwise. To begin the proof of \netheo{main_light} we derive the LST of the stationary workload $\vk Q$ for any fixed output rates $\vk r>0$. The following result, whose proof is given in Appendix \ref{A3}, generalizes \cite[Theorem 6.1]{DDR} for queuing networks. Define,
for any $j\in\sprod{n-1}$, $\vk\omega\in\R_+^{n}$ and $s\in\R_+$,
\begin{align}
    \delta_j(\vk\omega) &:=  \sum_{l\in{\tt S}_j}\frac{\hat{p}_{l}\omega_l}{\hat{p}_j},\label{delta_form}\\
    \hat{\delta}_j(\vk\omega) &:= 
    \sum_{l\in{\tt S}_{j+1}}\frac{\hat{p}_{l}\omega_l}{\hat{p}_j},\label{hat_delta_form}\\
    \kappa_{j+1}(\vk\omega) &:=  \sum_{l=j+1}^{n}\left(\frac{r_{l-1}}{\hat{p}_{l-1}} - \frac{r_{l}}{\hat{p}_{l}}\right)\sum_{i\in{\tt S}_{l}}\hat{p}_{i}\omega_i,\label{kappa_form_1}\\
    \psi_i(s) &:= \log\E{e^{-sX_i(1)}} = r_is + \log\E{e^{-s\hat{p}_jJ(1)}},\label{psi_def}
\end{align}
and $\Phi_i(\cdot)$ is the inverse of $\psi_i(\cdot)$, Moreover,
\begin{align}
    \kappa_{j+1}(\vk\omega) & = \sum_{i=j+1}^{n}\left(\frac{r_{\max(j,i^{\prime})}}{\hat{p}_{\max(j,i^{\prime})}}-\frac{r_i}{\hat{p_i}}\right)\hat{p}_i\omega_i.\label{kappa_form_2}
\end{align}

\begin{lem}\label{lem:Laplace_for_networks} Consider the queuing network $(J,P,\vk r)$, and grant Assumptions $\mbox{\textbf{\bf P}}$, $\mbox{\textbf{\bf N1}}$, {$\mbox{\textbf{\bf N2}}$ and $\mbox{\textbf{\bf R}}$}. Then, for any $\vk \omega\in\R_+^n$
    \begin{align}
    \E{e^{-\sprod{\vk\omega,\vk Q}}} =-\E{X_n(1)}\frac{\omega_n}{\psi_n(\omega_n)}\prod_{j=1}^{n-1}\frac{\Phi_j\bigl(\kappa_{j+1}(\vk\omega)\bigr) - \delta_j(\vk\omega)}{\Phi_j\bigl(\kappa_{j+1}(\vk\omega)\bigr) - \hat{\delta}_j(\vk\omega)}\frac{\kappa_{j+1}(\vk\omega) - \psi_j\bigl(\hat{\delta}_j(\vk\omega)\bigr)}{\kappa_{j+1}(\vk\omega) - \psi_j\bigl(\delta_j(\vk\omega)\bigr)}.\label{Laplace_W_simplify}
\end{align}
\end{lem}

Define, for a given $\vk\omega$,
\begin{align*}
    N_{j,1}(u) &:= \Phi_j\left(\kappa^{\star}_{j+1}(u)\right) - \delta^{\star}_j(u), \quad N_{j,2}(u):= \kappa^{\star}_{j+1}(u) - \psi_j\bigl(\hat{\delta}^{\star}_j(u)\bigr),\\
    D_{j,1}(u) &:=\Phi_j\left(\kappa^{\star}_{j+1}(u)\right) - \hat{\delta}^{\star}_j(u),\quad D_{j,2}(u):=\kappa^{\star}_{j+1}(u) - \psi_j\bigl(\delta^{\star}_j(u)\bigr),
\end{align*}
where we denote, for any $j\in\sprod{n-1}$, 
\begin{align*}
    \kappa^{\star}_{j+1}(u):= \kappa_{j+1}\bigl(\vk\omega^{\star}(u)\bigr),\qquad \delta^{\star}_{j}(u) := \delta_{j}\bigl(\vk\omega^{\star}(u)\bigr),\qquad \hat{\delta}^{\star}_{j}(u) := \hat{\delta}_j\bigl(\vk\omega^{\star}(u)\bigr),
\end{align*}
with, any $j\in\sprod{n}$,
\begin{align}
    \omega^{\star}_{j}(u) := r^{\beta}_j(u)\,\omega_j.\label{omega_star_def}
\end{align}
As a direct consequence of the above lemma, by \eqref{Laplace_W_simplify}, the joint LST of the random vector $\vk Q^{r,u}$ equals
\begin{align}
    \E{e^{-\sprod{\vk\omega,\vk Q^{r,u}}}} &=\E{e^{-\sprod{\vk\omega^{\star}(u), \vk{Q}_u}}}=-\E{X_n(1)}\frac{\omega^{\star}_n(u)}{\psi_n\bigl(\omega^{\star}_n(u)\bigr)}\prod_{j=1}^{n-1}\frac{N_{j,1}(u)}{D_{j,1}(u)}\frac{N_{j,2}(u)}{D_{j,2}(u)}\label{Laplace_W_in_r_c}.
\end{align}
Using \eqref{delta_form}, \eqref{hat_delta_form} \eqref{kappa_form_1} and \eqref{kappa_form_2} it follows that, for any $j\in\sprod{n-1}$,
\begin{align}
    \delta^{\star}_{j}(u) &= \sum_{l\in{\tt S}_j}\frac{\hat{p}_{l}\omega_lr_{l}^{\beta}(u)}{\hat{p}_j},\label{delta_r_def}\\
    \hat{\delta}^{\star}_{j}(u) &= \sum_{l\in{\tt S}_{j+1}}\frac{\hat{p}_{l}\omega_l r_{l}^{\beta}(u)}{\hat{p}_j},\label{delta_r_hat_def}\\
    \kappa^{\star}_{j+1}(u) &= \sum_{l=j+1}^{n}\left(\frac{r_{l-1}(u)}{\hat{p}_{l-1}} - \frac{r_l(u)}{\hat{p}_{l}}\right)\sum_{i\in{\tt S}_{l}}\hat{p}_{i}\omega_i r_{i}^{\beta}(u),\label{kappa_r_def}
\end{align}

We separately assess the numerator and denominator of each of the factors appearing in the product in \eqref{Laplace_W_in_r_c} (in the regime that $u\to\infty$), but we start with the denominator of the fraction in front of the product.

\medskip

\noindent {\it --- Expansion of $\psi_n\bigl(\omega_n^{\star}(u)\bigr)$.}
Assumption \textbf{L} combined with the definition \eqref{psi_def} tells us that, as $s\to\infty$, uniformly in $u$,
\begin{align}
    \psi_j(s) &= \log\E{e^{-sX_j(1)}} = \log\E{e^{-s(\hat{p}_j J(1) - r_j(u))}} = \log\E{e^{sr_j(u))}}\log\E{e^{-s\hat{p}_j J(1)}}\notag \\
    &= r_j(u)s + \mathfrak{C}\hat{p}_j^{\alpha}s^{\alpha} + o(s^{\alpha}).\label{psi_repr_c}
\end{align}
Hence, recalling \eqref{omega_star_def} and the fact that $\beta+1 = \alpha\beta$,
\begin{align}
    \psi_n\bigl(\omega_n^{\star}(u)\bigr) &= r_n(u)r_n^{\beta}(u)\omega_n + \mathfrak{C}\hat{p}_n^{\alpha}\left(r_n^{\beta}(u)\omega_n\right)^{\alpha} + o\left(r_n^{\beta}(u)\right)\notag\\
    &=\left(\omega_n + \mathfrak{C}\hat{p}_n^{\alpha}\omega_n^{\alpha}\right)r_n^{\alpha\beta}(u)+ o\left(r_n^{\alpha\beta}(u)\right)\label{psi_omega_star_exp_c}
\end{align}

--- {\it Expansion of $D_{j,2}(u)$.} Applying \eqref{r_comp} to \eqref{delta_r_def} and \eqref{delta_r_hat_def}, we obtain, for any $j\in\sprod{n-1}$,
\begin{align}
    \delta_j^{\star}(u) &= \mathcal{r}^{\beta}_{\mathcal{k}(j)}(u)\sum_{l\in{\tt S}^{\star}_j}\frac{\hat{p}_{l}\omega_l\mathfrak{r}_{l}^{\beta}}{\hat{p}_j} +o\left(\mathcal{r}^{\beta}_{\mathcal{k}(j)}(u)\right),\label{delta_exp_c}\\
    \hat{\delta}_j^{\star}(u) &= \mathcal{r}^{\beta}_{\mathcal{k}(j+1)}(u)\sum_{l\in{\tt S}^{\star}_{j+1}}\frac{\hat{p}_{l}\omega_l\mathfrak{r}_{l}^{\beta}}{\hat{p}_j} +o\left(\mathcal{r}^{\beta}_{\mathcal{k}(j+1)}(u)\right),\label{delta_hat_exp_c}
\end{align}
and, combining \eqref{kappa_r_def} with \eqref{r_comp} we have, for any $j\in\sprod{n-1}$,
\begin{align}
    \kappa_{j+1}^{\star}(u) = f_j\mathcal{r}_{\mathcal{k}(j)}(u)\mathcal{r}^{\beta}_{\mathcal{k}(j+1)}(u)+o\left(\mathcal{r}_{\mathcal{k}(j)}(u)\mathcal{r}^{\beta}_{\mathcal{k}(j+1)}(u)\right),\label{kappa_exp_2_c}
\end{align}
where, recalling the definition \eqref{F_cal},
\begin{align} 
f_j &= \sum_{l=j+1}^{n}\left(\mathbb{I}_{\{l-1\in\mathcal{k}(j)\}}\frac{\mathfrak{r}_{l-1}}{\hat{p}_{l-1}} - \mathbb{I}_{\{l\in\mathcal{k}(j+1)\}}\frac{\mathfrak{r}_l}{\hat{p}_{l}}\right)\sum_{i\in{\tt S}^{\star}_{l}}\hat{p}_{i}\omega_i \mathfrak{r}_{i}^{\beta}\notag\\
&=\begin{cases}
    \sum\limits_{\substack{l=j+1 \\ l\in{\tt I}_{\mathcal{k}(j+1)}}}^{n}\left(\frac{\mathfrak{r}_{l-1}}{\hat{p}_{l-1}} - \frac{\mathfrak{r}_l}{\hat{p}_{l}}\right)\sum\limits_{i\in{\tt S}^{\star}_{l}}\hat{p}_{i}\omega_i\mathfrak{r}_i^{\beta},\qquad &\text{ if } \mathcal{k}(j) = \mathcal{k}(j+1),\\
    \frac{\mathfrak{r}_{j}}{\hat{p}_j}\sum\limits_{l\in{\tt S}_{j+1}^{\star}}\hat{p}_l\omega_l\mathfrak{r}_l^{\beta}, \qquad & \text{ if } \mathcal{k}(j) \not= \mathcal{k}(j+1),
\end{cases}\notag\\
&=\begin{cases}
\mathcal{F}_{\alpha,j}\left(\vk{\mathfrak{r}}_{{\tt I}_k}\circ\vk\omega_{{\tt I}_{\mathcal{k}(j)}}\right), \qquad &\text{ if } \mathcal{k}(j) = \mathcal{k}(j+1),\\
\frac{\mathfrak{r}_{j}}{\hat{p}_j}\sum\limits_{l\in{\tt S}_{j+1}^{\star}}\hat{p}_l\omega_l\mathfrak{r}_l^{\beta}, \qquad\qquad &\text{ if } \mathcal{k}(j) \not= \mathcal{k}(j+1).
\end{cases}\notag
\end{align}

    Upon combining \eqref{psi_repr_c} with \eqref{delta_exp_c}, we obtain, for any $j\in\sprod{n-1}$,
\begin{align}
    \psi_j\left(\delta_j^{\star}(u)\right) &= r_j(u)\delta_j^{\star}(u) + \mathfrak{C}\hat{p}_j^{\alpha} \bigl(\delta_j^{\star}(u)\bigr)^{\alpha} + o\left(\bigl(\delta_j^{\star}(u)\bigr)^{\alpha}\right)\notag\\
    &= \left(\mathfrak{r}_{j}\left(\sum_{l\in{\tt S}^{\star}_j}\frac{\hat{p}_{l}\omega_l\mathfrak{r}_{l}^{\beta}}{\hat{p}_j}\right) + \mathfrak{C}\hat{p}_j^{\alpha} \left(\sum_{l\in{\tt S}^{\star}_j}\frac{\hat{p}_{l}\omega_l\mathfrak{r}_{l}^{\beta}}{\hat{p}_j}\right)^{\alpha} \right)\mathcal{r}_{\mathcal{k}(j)}^{\alpha\beta}(u) + o\left(\mathcal{r}_{\mathcal{k}(j)}^{\alpha\beta}(u)\right).\label{psi_delta_exp_c}
\end{align}
Using \eqref{S_rec} we know that, for any $j\in\sprod{n-1}$,
\begin{align}
    \sum\limits_{l\in{\tt S}_{j+1}^{\star}}\hat{p}_l\omega_l\mathfrak{r}_l^{\beta} - \sum\limits_{l\in{\tt S}^{\star}_j}\hat{p}_l\omega_l\mathfrak{r}_l^{\beta} = \sum\limits_{l\in{\tt D}_{j}^{\star}}\hat{p}_l\omega_l\mathfrak{r}_l^{\beta} -\hat{p}_j\omega_j\mathfrak{r}_j^{\beta}.\label{prep_calc}
\end{align}

Hence, by combining \eqref{kappa_exp_2_c} with \eqref{psi_delta_exp_c}, we conclude that  the following expansion is valid: as $u\to\infty$, for any $j\in\sprod{n-1}$, 
\begin{align}
    D_{j,2}(u)=\kappa^{\star}_{j+1}(u) - \psi_{j}\bigl(\delta_j^\star(u)\bigr)&= a_{j}\mathcal{r}_{\mathcal{k}(j)}^{\alpha\beta}(u) + o\left(\mathcal{r}_{\mathcal{k}(j)}^{\alpha\beta}(u)\right),\label{kappa-psi_delta_exp_c}
\end{align}
where, due to \eqref{prep_calc}, for any $j\in\sprod{n-1}$,
\begin{align}
    a_{j} &:= \mathbb{I}_{\mathcal{k}(j) = \mathcal{k}(j+1)}f_j - \mathfrak{r}_{j}\sum_{l\in{\tt S}^{\star}_j}\frac{\hat{p}_{l}\omega_l\mathfrak{r}_{l}^{\beta}}{\hat{p}_j} - \mathfrak{C} \left(\sum_{l\in{\tt S}^{\star}_j}\hat{p}_{l}\omega_l\mathfrak{r}_{l}^{\beta}\right)^{\alpha}\notag\\
    &=\begin{cases}
     \sum\limits_{\substack{l=j+1 \\ l\in{\tt I}_{\mathcal{k}(j+1)}}}^{n}\left(\frac{\mathfrak{r}_{l-1}}{\hat{p}_{l-1}} - \frac{\mathfrak{r}_l}{\hat{p}_{l}}\right)\sum\limits_{i\in{\tt S}^{\star}_{l}}\hat{p}_{i}\omega_i\mathfrak{r}_i^{\beta} - \mathfrak{r}_{j}\sum\limits_{l\in{\tt S}^{\star}_j}\frac{\hat{p}_{l}\omega_l\mathfrak{r}_{l}^{\beta}}{\hat{p}_j} - \mathfrak{C} \left(\sum\limits_{l\in{\tt S}^{\star}_j}\hat{p}_{l}\omega_l\mathfrak{r}_{l}^{\beta}\right)^{\alpha} ,  &\text{ if } \mathcal{k}(j) = \mathcal{k}(j+1),\\
    -\omega_j\mathfrak{r}_{j}^{\alpha\beta} - \mathfrak{C}\hat{p}_j^{\alpha} \omega^{\alpha}_j\mathfrak{r}_{j}^{\alpha\beta}, &\text{ if } \mathcal{k}(j) \not= \mathcal{k}(j+1).
    \end{cases}\notag\\
     &=\begin{cases}
     \sum\limits_{\substack{l=j+1 \\ l\in{\tt I}_{\mathcal{k}(j+1)}}}^{n-1}\frac{\mathfrak{r}_{l}}{\hat{p}_{l}} \left(\sum\limits_{i\in{\tt S}^{\star}_{l+1}}\hat{p}_{i}\omega_i\mathfrak{r}_i^{\beta} - \sum\limits_{i\in{\tt S}^{\star}_{l}}\hat{p}_{i}\omega_i\mathfrak{r}_i^{\beta}\right) - \omega_n\mathfrak{r}_n^{\alpha\beta} - \mathfrak{C} \left(\sum\limits_{l\in{\tt S}^{\star}_j}\hat{p}_{l}\omega_l\mathfrak{r}_{l}^{\beta}\right)^{\alpha} ,  &\text{ if } \mathcal{k}(j) = \mathcal{k}(j+1),\\
    -\omega_j\mathfrak{r}_{j}^{\alpha\beta} - \mathfrak{C}\hat{p}_j^{\alpha} \omega^{\alpha}_j\mathfrak{r}_{j}^{\alpha\beta}, &\text{ if } \mathcal{k}(j) \not= \mathcal{k}(j+1).
    \end{cases}\notag\\
    &=\begin{cases}
     \sum\limits_{\substack{l=j \\ l\in{\tt I}_{\mathcal{k}(j)}}}^{n-1}\frac{\mathfrak{r}_{l}}{\hat{p}_{l}}\sum\limits_{m\in{\tt D}^{\star}_l}\hat{p}_m\omega_m\mathfrak{r}_m^{\beta} - \sum\limits_{\substack{l=j \\ l\in{\tt I}_{\mathcal{k}(j)}}}^{n}\omega_l\mathfrak{r}_l^{\alpha\beta} - \mathfrak{C} \left(\sum\limits_{l\in{\tt S}^{\star}_j}\hat{p}_{l}\omega_l\mathfrak{r}_{l}^{\beta}\right)^{\alpha}, \qquad &\text{ if } \mathcal{k}(j) = \mathcal{k}(j+1),\\
    -\omega_j\mathfrak{r}_{j}^{\alpha\beta} - \mathfrak{C}\hat{p}_j^{\alpha} \omega^{\alpha}_j\mathfrak{r}_{j}^{\alpha\beta},  \qquad &\text{ if } \mathcal{k}(j) \not= \mathcal{k}(j+1).
    \end{cases}\notag\\
    &= \sum\limits_{\substack{l=j \\ l\in{\tt I}_{\mathcal{k}(j)}}}^{n}\frac{\mathfrak{r}_{l}}{\hat{p}_{l}}\sum\limits_{m\in{\tt D}^{\star}_l}\hat{p}_m\omega_m\mathfrak{r}_m^{\beta} - \sum\limits_{\substack{l=j \\ l\in{\tt I}_{\mathcal{k}(j)}}}^{n}\omega_l\mathfrak{r}_l^{\alpha\beta} - \mathfrak{C} \left(\sum\limits_{l\in{\tt S}^{\star}_j}\hat{p}_{l}\omega_l\mathfrak{r}_{l}^{\beta}\right)^{\alpha},\label{A_def_1}
\end{align}
as ${\tt D}^{\star}_j = \varnothing$ if $\mathcal{k}(j)\not= \mathcal{k}(j+1)$.

\begin{remark}
It is noted that the formula \eqref{A_def_1} is well-defined for $j=n$ as well. Hence, in the same way we can define the constant $a_n$.
\end{remark}
\COM{\begin{remark}
$A_j = A_j\bigl(\vk\omega_{{\tt I}_{\mathcal{k}(j)}}\bigr)$. Moreover, if ${\tt D}_l^{\star}=\varnothing$ for any $l\in{\tt I}_{\mathcal{k}(j)}$, then $A_j<0$ for any $\vk\omega\in(0,\infty)^{n}$. Otherwise, $A_j\not= 0$ almost everywhere because $A_j$ is linear with respect to $\omega_m$ for any $m\in\bigcup\limits_{\substack{l=j \\ l\in{\tt I}_{\mathcal{k}(j)}}}^{n}{\tt D}_{l}^{\star}$.
\end{remark}}

{\it --- Expansion of $N_{j,2}(u)$.}
First, combining \eqref{delta_r_hat_def} with \eqref{kappa_r_def} and using \eqref{prep_calc} we obtain, for any $j\in\sprod{n-1}$, recalling that ${\tt D}_n = \varnothing$,
\begin{align}
    \kappa^{\star}_{j+1}(u) - r_j(u)\hat{\delta}^{\star}_j(u) &= \sum_{l=j+1}^{n}\left(\frac{r_{l-1}(u)}{\hat{p}_{l-1}} - \frac{r_l(u)}{\hat{p}_{l}}\right)\sum_{i\in{\tt S}_{l}}\hat{p}_{i}\omega_i r_{i}^{\beta}(u) - \frac{r_j(u)}{\hat{p}_j}\sum_{l\in{\tt S}_{j+1}}\hat{p}_{l}\omega_l r_{l}^{\beta}(u)\notag\\
    &=\sum_{l=j+2}^{n}\frac{r_{l-1}(u)}{\hat{p}_{l-1}} \sum_{i\in{\tt S}_{l}}\hat{p}_{i}\omega_i r_{i}^{\beta}(u) - \sum_{l=j+1}^{n}\frac{r_l(u)}{\hat{p}_{l}}\sum_{i\in{\tt S}_{l}}\hat{p}_{i}\omega_i r_{i}^{\beta}(u)\notag\\
    &=\sum_{l=j+1}^{n-1}\frac{r_{l}(u)}{\hat{p}_{l}} \left(\sum_{i\in{\tt S}_{l+1}}\hat{p}_{i}\omega_i r_{i}^{\beta}(u) - \sum_{i\in{\tt S}_{l}}\hat{p}_{i}\omega_i r_{i}^{\beta}(u)\right) - \omega_n r_{n}^{\alpha\beta}(u)\notag\\
    &=\sum_{l=j+1}^{n-1}\frac{r_{l}(u)}{\hat{p}_{l}} \left(\sum_{m\in{\tt D}_{l}}\hat{p}_{m}\omega_m r_{m}^{\beta}(u) - \hat{p}_{l}\omega_l r_{l}^{\beta}(u)\right) - \omega_n r_{n}^{\alpha\beta}(u)\notag\\
    &=\sum_{l=j+1}^{n}\frac{r_{l}(u)}{\hat{p}_{l}} \sum_{m\in{\tt D}_{l}}\hat{p}_{m}\omega_m r_{m}^{\beta}(u) - \sum_{l=j+1}^{n}\omega_l r_{l}^{\alpha\beta}(u)\label{kappa-delta}.
\end{align}

Hence, applying \eqref{psi_repr_c} and \eqref{kappa-delta} we obtain for any $j\in\sprod{n-1}$ the asymptotics of the denominator of the last fraction of \eqref{Laplace_W_in_r_c}:
\begin{align}
    N_{j,2}(u)=\kappa^{\star}_{j+1}(u) - \psi_{j}\bigl(\hat{\delta}_j^{\star}(u)\bigr) &= \kappa^{\star}_{j+1}(u) - r_j(u)\hat{\delta}_j^{\star}(u)  - \mathfrak{C}\hat{p}_j^{\alpha}\bigl(\hat{\delta}_j^{\star}(u)\bigr)^{\alpha} + o\bigl(\bigl(\hat{\delta}_j^{\star}(u)\bigr)^{\alpha}\bigr) \notag\\
    &= \sum_{l=j+1}^{n}\frac{r_{l}(u)}{\hat{p}_{l}} \sum_{i\in{\tt D}_{l}}\hat{p}_{i}\omega_i r_{i}^{\beta}(u) - \sum_{l=j+1}^{n}\omega_l r_{l}^{\alpha\beta}(u) - \mathfrak{C}\hat{p}_j^{\alpha}\bigl(\hat{\delta}_j^{\star}(u)\bigr)^{\alpha} + o\bigl(\bigl(\hat{\delta}_j^{\star}(u)\bigr)^{\alpha}\bigr)\notag\\
    &=b_j\mathcal{r}_{\mathcal{k}(j+1)}^{\alpha\beta} + o\bigl(\mathcal{r}_{\mathcal{k}(j+1)}^{\alpha\beta}\bigr),\label{kappa-psi_hat_delta_exp_c}
\end{align}
where, using \eqref{r_comp} and \eqref{delta_hat_exp_c}, 
\begin{align*}
    b_j &:= \sum_{\substack{l=j+1 \\ l\in{\tt I}_{\mathcal{k}(j+1)}}}^{n}\frac{\mathfrak{r}_{l}}{\hat{p}_{l}} \sum_{m\in{\tt D}^{\star}_{l}}\hat{p}_{m}\omega_m \mathfrak{r}_{m}^{\beta} - \sum_{\substack{l=j+1 \\ l\in{\tt I}_{\mathcal{k}(j+1)}}}^{n}\omega_l \mathfrak{r}_{l}^{\alpha\beta} - \mathfrak{C}\hat{p}_j^{\alpha}\left(\sum_{l\in{\tt S}^{\star}_{j+1}}\frac{\hat{p}_{l}\omega_l\mathfrak{r}_{l}^{\beta}}{\hat{p}_j}\right)^{\alpha}= a_{j+1}
\end{align*}
due to \eqref{A_def_1}.

\medskip 

{\it --- Expansion of $N_{j,1}(u)$ and $D_{j,1}(u)$.} Define for any $j\in\sprod{n-1}$
\begin{align}
    \phi_j(u) := \Phi_j\bigl(\kappa^{\star}_{j+1}(u)\bigr)\label{phi_rep_1_c}.
\end{align}
We need to derive the asymptotics of function $\phi_j(u)$ as $u\to\infty$ for every $j\in\sprod{n-1}$. On the one hand, using \eqref{kappa_exp_2_c},
\begin{align}
    \psi_j\bigl(\phi_j(u)\bigr) &= \kappa^{\star}_{j+1}(u) = f_{j}\mathcal{r}_{\mathcal{k}(j)}(u)\mathcal{r}_{\mathcal{k}(j+1)}^{\beta}(u) + o\left(\mathcal{r}_{\mathcal{k}(j)}(u)\mathcal{r}_{\mathcal{k}(j+1)}^{\beta}(u)\right),\label{psi_phi_1_c}
\end{align}
and on the other hand, using \eqref{psi_repr_c}
\begin{align}
    \psi_j\bigl(\phi_j(u)\bigr) = r_j(u)\phi_j(u) + \mathfrak{C}\hat{p}_j^{\alpha}\bigl(\phi_j(u)\bigr)^{\alpha} +  o\left(\bigl(\phi_j(u)\bigr)^{\alpha}\right).\label{psi_phi_2_c}
\end{align}

Comparing the leading terms of \eqref{psi_phi_1_c} and \eqref{psi_phi_2_c} we obtain that there may be only two cases: either
\begin{align*}
    \mathcal{r}_{\mathcal{k}(j)}(u)\mathcal{r}_{\mathcal{k}(j+1)}^{\beta}(u) \asymp r_j(u)\phi_j(u), 
\end{align*}
or
\begin{align*}
    \mathcal{r}_{\mathcal{k}(j)}(u)\mathcal{r}_{\mathcal{k}(j+1)}^{\beta}(u) \asymp \bigl(\phi_j(u)\bigr)^{\alpha}.
\end{align*}

The first case is equivalent to
\begin{align*}
    \phi_j(u)\asymp \mathcal{r}_{\mathcal{k}(j+1)}^{\beta}(u),
\end{align*}
and the second case is equivalent to
\begin{align*}
    \phi_j(u)\asymp \mathcal{r}_{\mathcal{k}(j)}^{\frac{1}{\alpha}}(u)\mathcal{r}_{\mathcal{k}(j+1)}^{\frac{\beta}{\alpha}}(u).
\end{align*}
If $\mathcal{k}(j) = \mathcal{k}(j+1)$, then this two cases are equivalent. However, if $\mathcal{k}(j) \not= \mathcal{k}(j+1)$ in the second case we can see that, applying \eqref{r_diff},
\begin{align*}
    \kappa_{j+i}^{\star}(u) = o\bigl(r_j(u)\phi_j(u)\bigr),
\end{align*}
implying that the second case is impossible. Thus, without loss of generality,
\begin{align}
    \phi_j(u) = g_{j}\mathcal{r}_{\mathcal{k}(j+1)}^{\beta}(u) + o\left(\mathcal{r}_{\mathcal{k}(j+1)}^{\beta}(u)\right),\label{phi_exp_c}
\end{align}
where $g_{j}$ is the solution of
\begin{align*}
 f_j = \mathfrak{r}_jg_j + \mathbb{I}_{\mathcal{k}(j) = \mathcal{k}(j+1)} \mathfrak{C}\hat{p}_j^{\alpha}g_j^{\alpha}
\end{align*}

\begin{remark} As $f_j\geqslant 0$ for any $j\in\sprod{n-1}$, constant $g_j$ is well-defined and non-negative for all $j\in\sprod{n-1}$.  
\end{remark}

Applying \eqref{phi_exp_c} combined with \eqref{delta_exp_c} and \eqref{delta_hat_exp_c} we obtain, as $u\to\infty$,
\begin{align}
    N_{j,1}(u):=\Phi_j\bigl(\kappa^{\star}_{j+1}(u)\bigr) - \delta_j^{\star}(u) &= \phi_j(u) - \delta_j^{\star}(u) = c_j \mathcal{r}_{\mathcal{k}(j)}^{\beta}(u) +o\left(\mathcal{r}_{\mathcal{k}(j)}^{\beta}(u) \right),\label{Phi_kappa-delta_exp_c}\\
    D_{j,1}(u):=\Phi_j\bigl(\kappa^{\star}_{j+1}(u)\bigr) - \hat\delta_j^{\star}(u) &= \phi_j(u) - \hat\delta_j^{\star}(u) = d^{\star}_j \mathcal{r}_{\mathcal{k}(j+1)}^{\beta}(u) +o\left(\mathcal{r}_{\mathcal{k}(j+1)}^{\beta}(u)\right),\label{Phi_kappa-hat_delta_exp_a}
\end{align}
where
\begin{align*}
    c_j &:= \mathbb{I}_{\mathcal{k}(j) = \mathcal{k}(j+1)} g_j-\sum_{l\in{\tt S}^{\star}_j}\frac{\hat{p}_{l}\omega_l\mathfrak{r}_{l}^{\beta}}{\hat{p}_j},\quad \quad 
    d^{\star}_j := g_j -\sum_{l\in{\tt S}^{\star}_{j+1}}\frac{\hat{p}_{l}\omega_l\mathfrak{r}_{l}^{\beta}}{\hat{p}_j}.
\end{align*}

The asymptotic relations \eqref{Phi_kappa-delta_exp_c} and \eqref{Phi_kappa-hat_delta_exp_a} are meaningful only if constants $c_j$ and $d^{\star}_j$ are non-zero. One may check that $c_j\not=0$ and, if $\mathcal{k}(j) = \mathcal{k}(j+1)$, ${d_j^{\star}}\not=0$ {for almost all $\vk\omega\in\R_+^n$ with respect to the Lebesgue measure.} {Thus, the relation \eqref{Phi_kappa-delta_exp_c} and, in case $\mathcal{k}(j) = \mathcal{k}(j+1)$, the relation \eqref{Phi_kappa-hat_delta_exp_a} are applicable for almost all $\vk\omega\in\R_+^d$ and, similar to Remark \ref{rem:zero_D}, can be extended for all $\vk\omega\in\R_+^{d}$.} However, if $\mathcal{k}(j) \not= \mathcal{k}(j+1)$, we have, for any $\vk\omega\in\R_+^{n}$,
\begin{align*}
    g_j = \frac{f_j}{\mathfrak{r}_j} = \sum_{l\in{\tt S}^{\star}_{j+1}}\frac{\hat{p}_{l}\omega_l\mathfrak{r}_{l}^{\beta}}{\hat{p}_j},
\end{align*}
{implying that $d^{\star}_j = 0$. So the relation \eqref{Phi_kappa-hat_delta_exp_a} is reduced to $D_{j,1}(u) = o(\mathcal{r}_{\mathcal{k}(j+1)}^{\beta}(u))$. Thus, to derive the leading term for $D_{j,1}(u)$ in case $\mathcal{k}(j) \not= \mathcal{k}(j+1)$} we need to expand $\phi_j(u)$ further. Combining \eqref{phi_exp_c} with \eqref{delta_hat_exp_c}
\begin{align*}
    \phi_j(u) = \hat{\delta}_j^{\star}(u) + o\big(\mathcal{r}_{\mathcal{k}(j+1)}^{\beta}(u)\big).
\end{align*}
Denote $\varrho_j(u) := \phi_j(u) - \hat{\delta}_j^{\star}(u)$. Then, as $u\to\infty$, $\varrho_j(u) = o\left(\mathcal{r}_{\mathcal{k}(j+1)}^{\beta}(u)\right)$. According to \eqref{psi_repr_c},
\begin{align*}
  \psi_j\bigl(\phi_j(u)\bigr) = \kappa_{j+1}^{\star}(u) = r_{j}(u)\left(\hat{\delta}_j^{\star}(u) + \varrho_j(u)\right) + \mathfrak{C}\hat{p}_j^{\alpha}\left(\hat{\delta}_j^{\star}(u) + \varrho_j(u)\right)^{\alpha} + o\left(\left(\hat{\delta}_j^{\star}(u)\right)^{\alpha}\right).
\end{align*}

Taking into account \eqref{kappa-delta},
\begin{align*}
  \sum_{l=j+1}^{n}\frac{r_{l}(u)}{\hat{p}_{l}} \sum_{m\in{\tt D}_{l}}\hat{p}_{m}\omega_m r_{m}^{\beta}(u) - \sum_{l=j+1}^{n}\omega_l r_{l}^{\alpha\beta}(u)= r_{j}(u) \varrho_j(u) + \mathfrak{C}\hat{p}_j^{\alpha}\left(\hat{\delta}_j^{\star}(u) + \varrho_j(u)\right)^{\alpha} + o\left(\left(\hat{\delta}_j^{\star}(u)\right)^{\alpha}\right).
\end{align*}

Hence, using \eqref{r_comp}, as $u\to\infty$ (recall the definition \eqref{A_def_1})
\begin{align*}
    r_{j}(u)\varrho_j(u) &= \left(\sum_{\substack{l=j+1 \\ l\in{\tt I}_{\mathcal{k}(i+1)}}}^{n}\frac{\mathfrak{r}_{l}}{\hat{p}_{l}} \sum_{m\in{\tt D}^{\star}_{l}}\hat{p}_{m}\omega_m \mathfrak{r}_{m}^{\beta} - \sum_{\substack{l=j+1 \\ l\in{\tt I}_{\mathcal{k}(j+1)}}}^{n}\omega_l \mathfrak{r}_{l}^{\alpha\beta} - \mathfrak{C}\hat{p}_j^{\alpha}\left(\sum_{l\in{\tt S}^{\star}_{j+1}}\frac{\hat{p}_{l}\omega_l\mathfrak{r}_{l}^{\beta}}{\hat{p}_j}\right)^{\alpha}\right) \mathcal{r}^{\alpha\beta}_{\mathcal{k}(j+1)}(u)\\
    &\qquad\qquad+o\left(\mathcal{r}^{\alpha\beta}_{\mathcal{k}(j+1)}(u)\right)\\
    &= a_{j+1}\mathcal{r}^{\alpha\beta}_{\mathcal{k}(j+1)}(u)+o\left(\mathcal{r}^{\alpha\beta}_{\mathcal{k}(j+1)}(u)\right).
\end{align*}
 It implies that
\begin{align}
    \varrho_j(u) = d^{\star\star}_j\mathcal{r}_{\mathcal{k}(j+1)}^{\alpha\beta}(u)\mathcal{r}_{\mathcal{k}(j)}^{-1}(u) + o\left(\mathcal{r}_{\mathcal{k}(j+1)}^{\alpha\beta}(u)\mathcal{r}_{\mathcal{k}(j)}^{-1}(u)\right),\label{rho_exp}
\end{align}
for $d^{\star\star}_j := a_{j+1}/\mathfrak{r}_j$. Hence, in case $\mathcal{k}(j) \not= \mathcal{k}(j+1)$ we have that
\begin{align}
     {D_{j,1}(u)=\Phi_j(\kappa^{\star}_{j+1}) - \hat{\delta}_j^{\star}(u)} = \varrho_j(u) =  d^{\star\star}_j\mathcal{r}_{\mathcal{k}(j+1)}^{\alpha\beta}(u)\mathcal{r}_{\mathcal{k}(j)}^{-1}(u) + o\left(\mathcal{r}_{\mathcal{k}(j+1)}^{\alpha\beta}(u)\mathcal{r}_{\mathcal{k}(j)}^{-1}(u)\right).\label{Phi_kappa-hat_delta_exp_b}
\end{align}
We now combine the two asymptotic representations for $D_{j,1}(u)$: \eqref{Phi_kappa-hat_delta_exp_a} in case $\mathcal{k}(j) = \mathcal{k}(j+1)$  and \eqref{Phi_kappa-hat_delta_exp_b} in case $\mathcal{k}(j) \not= \mathcal{k}(j+1)$. We thus obtain the general version of the asymptotics as $u\to\infty$:
\begin{align}
    {D_{j,1}(u)} &= \,\bigr(d_j+o(1)\bigl)\mathcal{r}_{\mathcal{k}(j+1)}^{\beta}(u)\left(\frac{\mathcal{r}_{\mathcal{k}(j+1)}}{\mathcal{r}_{\mathcal{k}(j)}}\right)^{\mathbb{I}_{\mathcal{k}(j)\not= \mathcal{k}(j+1)}}\hspace{-2mm}=\, d_j\mathcal{r}_{\mathcal{k}(j+1)}^{\alpha\beta}(u)\mathcal{r}^{-1}_{\mathcal{k}(j)}+o\left(\mathcal{r}_{\mathcal{k}(j+1)}^{\alpha\beta}(u)\mathcal{r}^{-1}_{\mathcal{k}(j)}\right),\label{Phi_kappa-hat_delta_exp_c}
\end{align}
where
\begin{align*}
    d_j := \begin{cases}
        d^{\star}_j = g_j -\sum\limits_{l\in{\tt S}^{\star}_{j+1}}\frac{\hat{p}_{l}\omega_l\mathfrak{r}_{l}^{\beta}}{\hat{p}_j},\qquad\qquad &\text{ if } \mathcal{k}(j) = \mathcal{k}(j+1),\\
        d^{\star\star}_j = \frac{a_{j+1}}{\mathfrak{r}_j} \qquad\qquad &\text{ if } \mathcal{k}(j) \not= \mathcal{k}(j+1).
    \end{cases}
\end{align*}

{\it --- Combining the above expansions.} 
Inserting \eqref{psi_omega_star_exp_c}, \eqref{kappa-psi_delta_exp_c}, \eqref{kappa-psi_hat_delta_exp_c}, \eqref{Phi_kappa-delta_exp_c} and \eqref{Phi_kappa-hat_delta_exp_c} into \eqref{Laplace_W_in_r_c}, we obtain the final asymptotic representation for the LST of the vector $\vk Q^{r,u}$:
\begin{align*}
    \E{e^{-\sprod{\vk\omega,\vk Q^{r,u}}}} &= \E{e^{-\sprod{\vk\omega^{\star}(u),{\vk{Q}_u}}}}\\
    &= -\E{X_n(1)}\frac{\omega^{\star}_n(u)}{\psi_n\bigl(\omega^{\star}_n(u)\bigr)}\prod_{j=1}^{n-1}\frac{\Phi_j\left(\kappa^{\star}_{j+1}(u)\right) - \delta^{\star}_j(u)}{\Phi_j\left(\kappa^{\star}_{j+1}(u)\right) - \hat{\delta}^{\star}_j(u)}\frac{\kappa^{\star}_{j+1}(u) - \psi_j\bigl(\hat{\delta}^{\star}_j(u)\bigr)}{\kappa^{\star}_{j+1}(u) - \psi_j\bigl(\delta^{\star}_j(u)\bigr)}\\
    &\sim -r_n(u)\frac{r_n^{\beta}(u)\omega_n}{\left(\omega_n + \mathfrak{C}\hat{p}_n^{\alpha}\omega_n^{\alpha}\right)r_n^{\alpha\beta}(u)}\prod_{j=1}^{n-1}\frac{c_j \mathcal{r}_{\mathcal{k}(j)}^{\beta}(u)}{d_j\mathcal{r}_{\mathcal{k}(j+1)}^{\alpha\beta}(u)\mathcal{r}^{-1}_{\mathcal{k}(j)}}\frac{b_j\mathcal{r}_{\mathcal{k}(j+1)}^{\alpha\beta}}{a_{j}\mathcal{r}_{\mathcal{k}(j)}^{\alpha\beta}(u)}\\
    &=-\frac{\omega_n\mathfrak{r}_n^{\alpha\beta}}{\left(\omega_n + \mathfrak{C}\hat{p}_n^{\alpha}\omega_n^{\alpha}\right)\mathfrak{r}^{\alpha\beta}_n}\prod_{j=1}^{n-1}\frac{c_j}{a_j}\frac{b_j}{a_{j}}=-\frac{\omega_n\mathfrak{r}_n^{\alpha\beta} a_n}{\left(\omega_n + \mathfrak{C}\hat{p}_n^{\alpha}\omega_n^{\alpha}\right)\mathfrak{r}^{\alpha\beta}_n a_1}\prod_{j=1}^{n-1}\frac{c_j}{d_j}  \\&= \frac{\omega_n\mathfrak{r}_n^{\alpha\beta}}{a_1}\prod_{j=1}^{n-1}\frac{c_j}{d_j} = \prod_{k=1}^{m}f_k^{\star}
\end{align*}
where, for any $k\in\sprod{m}$,
\begin{align*}
    f^{\star}_k := \frac{\omega_{q_k+\abs{{\tt I}_k}-1}\mathfrak{r}^{\alpha\beta}_{q_k+\abs{{\tt I}_k}}-1}{\abs{a_{q_k+1}}}\prod_{j=0}^{\abs{{\tt I}_{k}}-2}\frac{\abs{c_{q_k+j}}}{\abs{d_{q_k+j}}}.
\end{align*}
Hence, the claim follows due to, for any $j\in\sprod{\abs{{\tt I}_k}-1}$,
\begin{align*}
    a_{q_k} &= \mathcal{A}_{\alpha,k}(\mathfrak{C},\vk{\mathfrak{r}}_{{\tt I}_k}\circ\vk\omega_{{\tt I}_k}),\quad 
    c_{q_k+j} = \mathcal{C}_{\alpha,q_k+j}(\mathfrak{C},\vk{\mathfrak{r}}_{{\tt I}_k}\circ\vk\omega_{{\tt I}_k}),\quad
    d_{q_k+j} = \mathcal{D}_{\alpha,q_k+j}(\mathfrak{C},\vk{\mathfrak{r}}_{{\tt I}_k}\circ\vk\omega_{{\tt I}_k}),
\end{align*}
using that in this case $\mathcal{k}(j) = \mathcal{k}(j+1)$. This completes the proof of \netheo{main_light}.
\QED

\appendix

\section{}\label{sec:appendix}
\subsection{\prooflem{S_D_relations}}\label{A1}
Claim \eqref{D_intersection} follows from the fact that merges are not allowed in our system; that is, if $p_{ij}>0$ and $p_{kl}>0$, then either $j\neq l$ or $i = k$. Hence, for any distinct $j,k\in\sprod{n}$, if $i\in{\tt D}_j\cap{\tt D}_k$ for some $i\in\sprod{n}$, then by the definition of the sets ${\tt D}_j$ and ${\tt D}_k$ we have $p_{ji}>0$ and $p_{ki}>0$. This contradicts the no-merging condition stated above. Therefore, \eqref{D_intersection} follows. 

We now proceed with establishing claim \eqref{S_D_intersection}. Assume it is not correct, i.e., there exist $i,j,l\in\sprod{n}$ such that $l\geqslant j$, $i\in{\tt S}_j$ and $i\in{\tt D}_l$. The last relation implies, by definition of the set ${\tt D}_l$, that $i^\prime = l\geqslant j$, which contradicts the relation $i\in{\tt S}_j$. Hence, \eqref{S_D_intersection} follows.

Consider next claim \eqref{S_D_union}. By definition, for any $j\in\sprod{n}$ and $l\in\sprod{n}^{j+}$ we have that ${\tt S}_j,\, {\tt D}_l\subset\sprod{n}^{j+}$, implying that ${\tt S}_{j}\cup\bigcup_{l=j}^{n}{\tt D}_{l} \subset \sprod{n}^{j+}$. On the other hand, for any $i\in\sprod{n}^{j+}$ there are two cases: either $i^\prime<j$ implying that $i\in{\tt S}_j$, or $i^\prime\geqslant j$ implying that $i\in {\tt D}_{i^\prime}\subset\bigcup_{l=j}^{n}{\tt D}_{l}$. In both cases we obtain $i\in{\tt S}_{j}\cup\bigcup_{l=j}^{n}{\tt D}_{l}$, hence $\sprod{n}^{j+}\subset {\tt S}_{j}\cup\bigcup_{l=j}^{n}{\tt D}_{l}$ implying \eqref{S_D_union}.

We then prove claim \eqref{S_rec}. Fix some $i\in{\tt S}_{j+1}$. By definition of the set ${\tt S}_{j+1}$, $i^\prime < j+1$. Hence, either $i^\prime<j$ implying that $i\in{\tt S}_j$, or $i^\prime = j$ implying that $i\in{\tt D}_j$. In both cases, $i\in {\tt S}_{j}\cup {\tt D}_{j}$. Using further that $j<j+1$, hence $j\not\in{\tt S}_{j+1}$, we have that ${\tt S}_{j+1} \subset {\tt S}_{j}\cup {\tt D}_{j}\setminus \{j\}$. On the other hand, if $i\in {\tt D}_{j}$ then $i>j$ and $i^\prime = j<j+1$, hence $i\in{\tt S}_{j+1}$. And if $i\in{\tt S}_j\setminus\{j\}$, then $i\geqslant j+1$ and $i^\prime<j<j+1$, Hence, again $i\in{\tt S}_{j+1}$. Combining the last two assertions we have that $ {\tt S}_{j}\cup {\tt D}_{j}\setminus \{j\}\subset{\tt S}_{j+1}$, hence \eqref{S_rec} follows.

Finally, we consider claim \eqref{k_in_S}. Let $k\in{\tt S}_j$. Then, by definition of set ${\tt S}_j$, $k\geqslant j$ and $k^\prime< j$, implying that $j\in[k^\prime+1,k]$. In the other direction, if $j\in[k^\prime+1,k]$, then $k\geqslant j$ and $k^\prime < j$, hence by definition it follows that $k\in{\tt S}_j$.
\QED

\subsection{Proof of Remark \ref{rem:zero_D}}\label{A2}
This remark is heavily based on the following lemma. {Its proof is elementary and is therefore omitted.}
\begin{lem}\label{frac_lim_lem} Let $g(\cdot),\,h(\cdot)$ be two continuous function $g,h:\R\to\R_+$ such that $\lim_{s\to\infty}g(s)\in\R_+$ and $h_\infty:=\lim_{s\to\infty}h(s)\in\R_+$. For any constants $A,B>0$, $\alpha>1$ we denote  $\varphi_{\alpha,A,B}(s) := As + Bs^{\alpha}$, for $s\geqslant0$. Then
\begin{align*}  
    \lim_{s\to\infty}\frac{\varphi_{\alpha,A,B}^{-1}(g(s)) - h(s)}{g(s) - \varphi_{\alpha,A,B}(h(s))} = \frac{1}{\varphi^{\prime}_{\alpha,A,B}\bigl(h_\infty\bigr)}\in(0,\infty).
\end{align*}
\end{lem}

Let $\omega(s):R\to\R_+^n$ such that $\lim_{s\to\infty}\omega(s)\in\R_+^n$, {and}  $A = \mathfrak{r}_j$, $B = \mathfrak{C}\hat{p}^{\alpha}_j$ (so that, by \eqref{Psi_def}, $\varphi_{\alpha,A,B}(\cdot) = \Psi_{\alpha,\mathfrak{C},j}(\cdot)$) and 
\begin{align*}
    &g(s) := \mathcal{F}_{q_k}(\vk\omega(s)),\quad 
    h(s) := \sum_{l\in{\tt S}^{\star}_{q_k}}\frac{\hat{p}_{l}\omega_l(s)}{\hat{p}_{q_k}}.
\end{align*}
It follows that, recalling the definition \eqref{C_cal},
\begin{align*}
    \Psi_{\alpha,\mathfrak{C},j}^{-1}(g(s)) - h(s) &= \mathcal{C}_{\alpha,q_k}(\mathfrak{C},\vk\omega(s)),
\end{align*}
and, {combining relation \eqref{k_in_S} with the definition} \eqref{A_cal},
\begin{align*}
    g(s) - \Psi_{\alpha,\mathfrak{C},j}(h(s)) &= \mathcal{F}_{q_k}(\vk\omega(s)) - \mathfrak{r}_j\sum_{l\in{\tt S}^{\star}_{q_k}}\frac{\hat{p}_{l}\omega_l(s)}{\hat{p}_j} - C\left(\sum_{l\in{\tt S}^{\star}_{q_k}}\hat{p}_{l}\omega_l(s)\right)^{\alpha}\\
    &= \sum_{\substack{l=q_k+1 \\ l\in{\tt I}_{\mathcal{k}(q_k+1)}}}^{n}\left(\frac{\mathfrak{r}_{l-1}}{\hat{p}_{l-1}} - \frac{\mathfrak{r}_l}{\hat{p}_{l}}\right)\sum_{i\in{\tt S}^{\star}_{l}}\hat{p}_{i}\omega_i(s) - \mathfrak{r}_j\sum_{l\in{\tt S}^{\star}_{q_k}}\frac{\hat{p}_{l}\omega_l(s)}{\hat{p}_j} - \mathfrak{C}\left(\sum_{l\in{\tt S}^{\star}_{q_k}}\hat{p}_{l}\omega_l(s)\right)^{\alpha}\\
    &= \sum_{l\in{\tt I}_{k}}^{n}\left(\frac{\mathfrak{r}_{l-1}}{\hat{p}_{l-1}} - \frac{\mathfrak{r}_l}{\hat{p}_{l}}\right)\sum_{i\in{\tt S}^{\star}_{l}}\hat{p}_{i}\omega_i(s) - \mathfrak{r}_j\sum_{l\in{\tt S}^{\star}_{q_k}}\frac{\hat{p}_{l}\omega_l(s)}{\hat{p}_j} - \mathfrak{C}\left(\sum_{l\in{\tt S}^{\star}_{q_k}}\hat{p}_{l}\omega_l(s)\right)^{\alpha}\\
    &=\mathcal{A}_{\alpha,k}(\mathfrak{C},\vk\omega(s)).
\end{align*}
Hence, applying \nelem{frac_lim_lem}, 
\begin{align*}
    \lim_{s\to\infty}\frac{\mathcal{C}_{\alpha,q_k}(\mathfrak{C},\vk\omega(s))}{\mathcal{A}_{\alpha,k}(\mathfrak{C},\vk\omega(s))}\in(0,\infty),
\end{align*}
thus establishing \eqref{rem_2}. We now proceed with proving the second claim. For any $k\in m$ and $j\in\sprod{\abs{{\tt I}_k}-3}$, by taking $A_1 = \mathfrak{r}_{q_k+j}$, $B_1 = \mathfrak{C}\hat{p}^{\alpha}_{q_k+j}$ (i.e., $ \varphi_{\alpha,A_1,B_1}(\cdot) = \Psi_{\alpha,\mathfrak{C},q_k+j}(\cdot)$)
\begin{align*}
    g_1(s) &:= \mathcal{F}_{q_k+j}\bigl(\omega(s)\bigr),\quad
    h_1(s) := \sum_{l\in{\tt S}^{\star}_{q_k+j+1}}\frac{\hat{p}_{l}\omega_l(s)}{\hat{p}_{q_k+j}},
\end{align*}
we obtain by applying \nelem{frac_lim_lem}, recalling the definition \eqref{D_cal},
\begin{align}
    \lim_{s\to\infty}\frac{\mathcal{D}_{\alpha,q_k+j}(\mathfrak{C},\vk\omega(s))}{\mathcal{E}_{k,j}(\vk\omega(s))}\in (0,\infty),\label{lim_1}
\end{align}
where
\begin{align*}
    \mathcal{E}_{k,j}(\vk\omega(s)) &= \mathcal{F}_{q_k+j}\bigl(\vk\omega(s)\bigr) - \mathfrak{r}_{q_k+j}h_1(s) - C\hat{p}_{q_k+j}^{\alpha}h_1^{\alpha}(s)\\
    &= \sum_{\substack{l=q_k+j+1 \\ l\in{\tt I}_{k}}}^{n}\left(\frac{\mathfrak{r}_{l-1}}{\hat{p}_{l-1}} - \frac{\mathfrak{r}_l}{\hat{p}_{l}}\right)\sum_{i\in{\tt S}^{\star}_{l}}\hat{p}_{i}\omega_i - \frac{\mathfrak{r}_{q_k+j}}{\hat{p}_{q_k+j}}\sum_{l\in{\tt S}^{\star}_{q_k+j+1}}\hat{p}_{l}\omega_l(s) - \mathfrak{C}\left(\sum_{l\in{\tt S}^{\star}_{q_k+j+1}}\hat{p}_{l}\omega_l(s)\right)^\alpha\\
    &=\sum_{\substack{l=q_k+j+2 \\ l\in{\tt I}_{k}}}^{n}\left(\frac{\mathfrak{r}_{l-1}}{\hat{p}_{l-1}} - \frac{\mathfrak{r}_l}{\hat{p}_{l}}\right)\sum_{i\in{\tt S}^{\star}_{l}}\hat{p}_{i}\omega_i - \frac{\mathfrak{r}_{q_k+j+1}}{\hat{p}_{q_k+j+1}}\sum_{l\in{\tt S}^{\star}_{q_k+j+1}}\hat{p}_{l}\omega_l(s) - \mathfrak{C}\left(\sum_{l\in{\tt S}^{\star}_{q_k+j+1}}\hat{p}_{l}\omega_l(s)\right)^\alpha.
\end{align*}
On the other hand, by taking $A_2=\mathfrak{r}_{q_k+j+1}$, $B_2 = \mathfrak{C}\hat{p}_{q_k+j+1}^{\alpha}$ (i.e., $\varphi_{\alpha,A_1,B_1}(\cdot) = \Psi_{\alpha,\mathfrak{C},q_k+j+1}(\cdot)$)
\begin{align*}
    g_2(s) &:=\mathcal{F}_{q_k+j+1}(\vk\omega(s)),\quad
    h_2(s) := \sum_{l\in{\tt S}^{\star}_{q_k+j+1}}\frac{\hat{p}_{l}\omega_l}{\hat{p}_{q_k+j+1}},
\end{align*}
and applying \nelem{frac_lim_lem}, now recalling the definition \eqref{C_cal},
\begin{align}
    \lim_{s\to\infty}\frac{\mathcal{C}_{\alpha,q_k+j+1}(\mathfrak{C},\vk\omega(s))}{\mathcal{E}^{\circ}_{k,j}(\vk\omega(s))}\in(0,\infty),\label{lim_2}
\end{align}
where
\begin{align*}
    \mathcal{E}^{\circ}_{k,j}(\vk\omega(s)) &= \mathcal{F}_{q_k+j+1}(\vk\omega(s)) - \mathfrak{r}_{q_k+j+1}h_2(s) - C\hat{p}_{q_k+j+1}^{\alpha}h_2^{\alpha}(s)\\
    &=\sum_{\substack{l=q_k+j+2 \\ l\in{\tt I}_{k}}}^{n}\left(\frac{\mathfrak{r}_{l-1}}{\hat{p}_{l-1}} - \frac{\mathfrak{r}_l}{\hat{p}_{l}}\right)\sum_{i\in{\tt S}^{\star}_{l}}\hat{p}_{i}\omega_i - \frac{\mathfrak{r}_{q_k+j+1}}{\hat{p}_{q_k+j+1}}\sum_{l\in{\tt S}^{\star}_{q_k+j+1}}\hat{p}_{l}\omega_l(s) - \mathfrak{C}\left(\sum_{l\in{\tt S}^{\star}_{q_k+j+1}}\hat{p}_{l}\omega_l(s)\right)^\alpha\\
    &=\mathcal{E}_{k,j}(\vk\omega(s)).
\end{align*}
Hence, for any $j\leqslant \abs{{\tt I}_k} - 3$, \eqref{rem_2} follows upon combining \eqref{lim_1} and \eqref{lim_2}. Finally, for $j=\abs{{\tt I}_k} - 2$ and any $\vk\omega\in\R_+^{n}$ (recall that ${\tt I}_k = \{q_k,\ldots,q_k +\abs{{\tt I}_k} - 1\}$ and ${\tt S}_l^{\star}\subset {\tt I}_k$ for any $l\in {\tt I}_k$),
\begin{align*}
    \mathcal{F}_{q_k + \abs{{\tt I}_k} - 2}(\vk\omega) = \sum_{\substack{l=q_k + \abs{{\tt I}_k} - 1 \\ l\in{\tt I}_{k}}}^{n}\left(\frac{\mathfrak{r}_{l-1}}{\hat{p}_{l-1}} - \frac{\mathfrak{r}_l}{\hat{p}_{l}}\right)\sum_{i\in{\tt S}^{\star}_{l}}\hat{p}_{i}\omega_i = \left(\frac{\mathfrak{r}_{q_k + \abs{{\tt I}_k} - 2}}{\hat{p}_{q_k + \abs{{\tt I}_k} - 2}} - \frac{\mathfrak{r}_{q_k + \abs{{\tt I}_k} - 1}}{\hat{p}_{q_k + \abs{{\tt I}_k} - 1}}\right)\hat{p}_{q_k + \abs{{\tt I}_k} - 1}\omega_{q_k + \abs{{\tt I}_k} - 1},
\end{align*}
and
\begin{align*}
    \Psi_{\alpha,\mathfrak{C},q_k+\abs{{\tt I}_k} - 2}\left(\sum_{l\in{\tt S}^{\star}_{q_k+\abs{{\tt I}_k} - 1}}\frac{\hat{p}_{l}\omega_l}{\hat{p}_{q_k+\abs{{\tt I}_k} - 2}}\right) = \frac{\mathfrak{r}_{q_k + \abs{{\tt I}_k} - 2}}{\hat{p}_{q_k + \abs{{\tt I}_k} - 2}}\hat{p}_{q_k + \abs{{\tt I}_k} - 1}\omega_{q_k + \abs{{\tt I}_k} - 1} + \mathfrak{C}\hat{p}^{\alpha}_{q_k+\abs{{\tt I}_k} - 1}\vk\omega_{q_k+\abs{{\tt I}_k} - 1}^{\alpha},
\end{align*}
implying that
\begin{align*}
     \Psi_{\alpha,\mathfrak{C},q_k+\abs{{\tt I}_k} - 2}\left(\sum_{l\in{\tt S}^{\star}_{q_k+\abs{{\tt I}_k} - 1}}\frac{\hat{p}_{l}\omega_l}{\hat{p}_{q_k+\abs{{\tt I}_k} - 2}}\right) > \frac{\mathfrak{r}_{q_k + \abs{{\tt I}_k} - 2}}{\hat{p}_{q_k + \abs{{\tt I}_k} - 2}}\hat{p}_{q_k + \abs{{\tt I}_k} - 1}\omega_{q_k + \abs{{\tt I}_k} - 1}  > \mathcal{F}_{q_k + \abs{{\tt I}_k} - 2}(\vk\omega).
\end{align*}
Hence, as $\Psi_{\alpha,\mathfrak{C},q_k+\abs{{\tt I}_k} - 2}(\cdot)$ is increasing,
\begin{align*}
    \mathcal{D}_{\alpha,\abs{{\tt I}_k} - 2}(C,\vk\omega) &=  \Psi_{\alpha,C,q_k+\abs{{\tt I}_k} - 2}^{-1}\bigl(\mathcal{F}_{q_k+\abs{{\tt I}_k} - 2}(\vk\omega)\bigr) -\sum_{l\in{\tt S}^{\star}_{\abs{{\tt I}_k} - 1}}\frac{\hat{p}_{l}\omega_l}{\hat{p}_j}<0,
\end{align*}
thus establishing the claim.
\QED

\subsection{\prooflem{lem:Laplace_for_networks}}\label{A3}
For any $j\in\sprod{n-1}$, $t\in\R$
\begin{align*}
    X_{j+1}(t) &= \hat{p}_{j+1}J(t) - r_{j+1}t = \hat{p}_{j+1}\frac{X_j(t) + r_j t}{\hat{p}_j} - r_{j+1}t \\
    &= \frac{\hat{p}_{j+1}}{\hat{p}_j}X_j(t) + \hat{p}_{j+1}\left(\frac{r_j}{\hat{p}_j} - \frac{r_{j+1}}{\hat{p}_{j+1}}\right)t\,=\, K_{j+1}X_j(t) + \Upsilon_{j+1}(t),
\end{align*}
where
\begin{align*}
    K_{j+1} := \frac{\hat{p}_{j+1}}{\hat{p}_{j}},\qquad \Upsilon_{j+1}(t) := \hat{p}_{j+1}\left(\frac{r_j}{\hat{p}_j} - \frac{r_{j+1}}{\hat{p}_{j+1}}\right)t.
\end{align*}

Note that $K_{j+1}>0$ and, due to Assumption {\textbf{N2}}, $\Upsilon_{j+1}(t)>0$. 
Recall that the process $J(\cdot)$ is centered and $r_j>0$, hence $X_j(t)\to-\infty$ as $t\to\infty$ almost surely, implying that Assumption \textbf{D} from \cite{DDR} is valid. Hence, applying \cite[Theorem 3.1]{DDR} for $\alpha = \vk 0$ and $\beta =\vk\omega$, the random vector $\overline{\vk X}$ has the following joint LST:
\begin{align}
    \E{e^{-\sprod{\vk\omega,\overline{\vk X}}}} = \E{e^{-\omega_n \overline{X}_n}}\prod_{j=1}^{n-1}\frac{\E{e^{-\left(\sum_{l=j+1}^{n}\theta_l^{\Upsilon}(\sum_{k=l}^{n}K_l^k\omega_k)\right)G_j - \left(\sum_{l=j}^{n}K_j^l\omega_{l}\right) \overline{X}_j}}}{\E{e^{-\left(\sum_{l=j+1}^{n}\theta_l^{\Upsilon}(\sum_{k=l}^{n}K_l^k\omega_k)\right)G_j - \left(\sum_{l=j+1}^{n}K_j^l\omega_{l}\right) \overline{X}_j}}}\label{X_bar_laplace},
\end{align}
where, using that $\Upsilon_{j}(t)$ is deterministic and linear,
\begin{align*}
    \theta_j^{\Upsilon}(t) &:= -\log\E{e^{-t\Upsilon_j(1)}} = t\Upsilon_{j}(1) = \Upsilon_{j}(t),\quad\:\:\:
    K_{i}^{j} := \prod_{l=i+1}^{j}K_l = \frac{\hat{p}_{j}}{\hat{p}_{i}}.
\end{align*}

and $\vk G \equiv (G_1,\ldots,G_n)$ with $G_j := \inf\{t\geqslant 0\colon X_j(t) =\overline{X}_j\text{ or }X_j(-t) =\overline{X}_j\}$ for any $j\in\sprod{n}$ is the first moment that the process $X(t)$ attains its maximum. Using \eqref{W_def},
\begin{align}
    \E{e^{-\sprod{\vk\omega,\vk Q}}} = \E{e^{-\sprod{\widetilde{\vk\omega}, \overline{\vk X}}}},\label{Laplace_W_to_X}
\end{align}
where, for any $j\in\sprod{n}$,
\begin{align*}
    \widetilde{\omega}_j:=\left((I-P)\vk\omega\right)_j =\omega_j - \sum_{i=1}^{n}p_{j,i}\omega_i = \omega_j - \sum_{i\in{\tt D}_j}p_{j,i}\omega_i.
\end{align*}
Combining \eqref{X_bar_laplace} with \eqref{Laplace_W_to_X} we obtain that
\begin{align}
    \E{e^{-\sprod{\vk\omega,\vk Q}}} = \E{e^{-\widetilde{\omega}_n \overline{X}_n}}\prod_{j=1}^{n-1}\frac{\E{e^{-\left(\sum_{l=j+1}^{n}\theta_l^{\Upsilon}(\sum_{k=l}^{n}K_l^k\widetilde{\omega}_k)\right)G_j - \left(\sum_{l=j}^{n}K_j^l\widetilde{\omega}_{l}\right) \overline{X}_j}}}{\E{e^{-\left(\sum_{l=j+1}^{n}\theta_l^{\Upsilon}(\sum_{k=l}^{n}K_l^k\tilde{\omega}_k)\right)G_j - \left(\sum_{l=j+1}^{n}K_j^l\widetilde{\omega}_{l}\right) \overline{X}_j}}}.\label{Laplace_W}
\end{align}
The next step is to simplify the equation \eqref{Laplace_W}. Denote
\begin{align*}
    \kappa^{\star}_{j+1}(\vk\omega) &:= \sum_{l=j+1}^{n}\theta_l^{\Upsilon}\left(\sum_{k=l}^{n}K_l^k\widetilde{\omega}_k\right), \quad
    \delta^{\star}_{j}(\vk\omega) := \sum_{l=j}^{n}K_j^l\widetilde{\omega}_{l},\quad
    \hat{\delta}^{\star}_{j}(\vk\omega):= \sum_{l=j+1}^{n}K_j^l\widetilde{\omega}_{l}.
\end{align*}
One can readily verify that, by using \eqref{S_D_strong_union}, 
\begin{align*}
    \delta^{\star}_j(\vk\omega) &=  \sum_{l=j}^{n}K_j^l\left(\omega_l - \sum_{i\in{\tt D}_l}p_{l,i}\omega_i\right) = \sum_{l=j}^{n}K_j^l\omega_l - \sum_{l=j}^{n}\sum_{i\in{\tt D}_l}K_{j}^{i}\omega_i = \sum_{l\in{\tt S}_j}K_j^l\omega_l = \delta_j(\vk\omega),\\
    \hat{\delta}^{\star}_j(\vk\omega) &= K_{j+1}\delta_{j+1} = \hat{\delta}_j(\vk\omega),\quad \quad\quad
    \kappa^{\star}_{j+1}(\vk\omega) = \sum_{l=j+1}^{n}\Upsilon_{l}(\delta_l) = \kappa_{j+1}(\vk\omega).
\end{align*}
Now using \eqref{k_in_S}, we can prove \eqref{kappa_form_2}:
\begin{align*}
    \kappa_{j+1}(\vk\omega) &= \sum_{l=j+1}^{n}\left(\frac{r_{l-1}}{\hat{p}_{l-1}} - \frac{r_{l}}{\hat{p}_{l}}\right)\sum_{i\in{\tt S}_{l}}\hat{p}_{i}\omega_i = \sum_{i=j+1}^{n}\hat{p}_{i}\omega_i\sum_{l=j+1}^{i}\mathbb{I}_{\{i\in{\tt S}_{l}\}}\left(\frac{r_{l-1}}{\hat{p}_{l-1}} - \frac{r_{l}}{\hat{p}_{l}}\right)\notag\\
    &=\sum_{i=j+1}^{n}\hat{p}_{i}\omega_i\sum_{l=j+1}^{i}\mathbb{I}_{\{ l\in[i^{\prime}+1,i]\}}\left(\frac{r_{l-1}}{\hat{p}_{l-1}} - \frac{r_{l}}{\hat{p}_{l}}\right) =\sum_{i=j+1}^{n}\hat{p}_{i}\omega_i\sum_{l=\max(j,i^{\prime})+1}^{i}\left(\frac{r_{l-1}}{\hat{p}_{l-1}} - \frac{r_{l}}{\hat{p}_{l}}\right)\notag\\
    & = \sum_{i=j+1}^{n}\left(\frac{r_{\max(j,i^{\prime})}}{\hat{p}_{\max(j,i^{\prime})}}-\frac{r_i}{\hat{p_i}}\right)\hat{p}_i\omega_i.
\end{align*}
Then apply \cite[Theorem VII.4]{bertoin1996levy}, which entails  that for any $x,y>0$
\begin{align*}
    \E{e^{-x G_i - y \overline{X}_i}} = -\E{X_i(1)}\frac{\Phi_i(x) - y}{x - \psi_i(y)},
\end{align*}
where $\psi_i(\cdot)$ is defined in \eqref{psi_def}, and $\Phi_i(\cdot)$ is the inverse of $\psi_i(\cdot)$. Hence, \eqref{Laplace_W_simplify} follows from \eqref{Laplace_W}.
\QED

\subsection{Proofs related to the example of a two-layer network}\label{A4} In this part of the appendix we prove \eqref{tree}. 
 In this case, for any $j\in\sprod{n}^{2+}$
\begin{align*}
    {\tt I}_1           &= \{1\}, 
    &\qquad 
    {\tt I}_2           &= \sprod{n}^{2+}, \\[6pt]
    \mathcal{r}_1(u)    &= r_1(u), 
    &\qquad 
    \mathcal{r}_2(u)    &= r_2(u), \\[6pt]
    \mathfrak{r}_1      &= 1,
    &\qquad 
    \hat{p}_1           &= 1,
    &\qquad 
    \hat{p}_{j}         &= p_j, \\[6pt]
    {\tt D}^{\star}_1   &= {\tt D}^{\star}_i = \varnothing,
    &\qquad 
    {\tt S}^{\star}_1   &= \{1\},
    &\qquad 
    {\tt S}^{\star}_j   &= \{j,\ldots,n\}.
\end{align*}
To apply \netheo{main_light} (or \netheo{main_heavy}, {respectively}) we first calculate the constants $\mathcal{A}_{\alpha,1}(\mathfrak{C},\omega_1)$ and $\mathcal{A}_{\alpha,2}(\mathfrak{C},\vk\omega)$ for $\vk\omega\in R_+^{\sprod{n}^{2+}}$ defined in \eqref{A_cal}
\begin{align*}
    \mathcal{A}_{\alpha,1}(\mathfrak{C},\omega_1) &= \sum\limits_{l\in{\tt I}_{1}}\frac{\mathfrak{r}_{l}}{\hat{p}_{l}}\sum\limits_{m\in{\tt D}^{\star}_l}\hat{p}_m\omega_m - \sum\limits_{l\in{\tt I}_{1}}\omega_l\mathfrak{r}_l - \mathfrak{C} \left(\sum\limits_{l\in{\tt S}^{\star}_{q_1}}\hat{p}_{l}\omega_l\right)^{\alpha} = -\omega_1 -\mathfrak{C}\omega_1^{\alpha},\\
    \mathcal{A}_{\alpha,2}(\mathfrak{C},\omega_1) &= \sum\limits_{l\in{\tt I}_{2}}\frac{\mathfrak{r}_{l}}{\hat{p}_{l}}\sum\limits_{m\in{\tt D}^{\star}_l}\hat{p}_m\omega_m - \sum\limits_{l\in{\tt I}_{2}}\omega_l\mathfrak{r}_l - \mathfrak{C} \left(\sum\limits_{l\in{\tt S}^{\star}_{q_2}}\hat{p}_{l}\omega_l\right)^{\alpha} = - \sum\limits_{l=2}^{n}\omega_l\mathfrak{r}_l - \mathfrak{C} \left(\sum\limits_{l=2}^{n}p_l\omega_l\right)^{\alpha}.
\end{align*}
Then we calculate the constants $\mathcal{F}_{\alpha,j}(\vk\omega)$ defined in \eqref{F_cal} for $j\in\sprod{n-1}^{2+}$:
\begin{align*}
    \mathcal{F}_{\alpha,j}(\omega) = \sum_{\substack{l=j+1 \\ l\in{\tt I}_{\mathcal{k}(j+1)}}}^{n}\left(\frac{\mathfrak{r}_{l-1}}{\hat{p}_{l-1}} - \frac{\mathfrak{r}_l}{\hat{p}_{l}}\right)\sum_{i\in{\tt S}^{\star}_{l}}\hat{p}_{i}\omega_i = \sum_{l=j+1 }^{n}\left(\frac{\mathfrak{r}_{l-1}}{p_{l-1}} - \frac{\mathfrak{r}_l}{p_{l}}\right)\sum_{i=l}^{n}p_{i}\omega_i = \sum_{l=j+1}^{n} \left(\frac{\mathfrak{r}_j}{p_j} - \frac{\mathfrak{r}_l}{p_l}\right)p_l\omega_l.
\end{align*}
Using further that for any $j\in\sprod{n-1}^{2+}$, according to \eqref{Psi_def},
\begin{align*}
    \Psi_{\alpha,\mathfrak{C},j}(s) = \mathfrak{r}_js + \mathfrak{C}\hat{p}_j^{\alpha}s^{\alpha} = \mathfrak{r}_js + \mathfrak{C}p_j^{\alpha}s^{\alpha} = \psi_j(s),
\end{align*}
we have for any $j\in\sprod{n-1}^{2+}$ (due to \eqref{C_cal} and \eqref{D_cal})
\begin{align*}
    \mathcal{C}_{\alpha,j}(\mathfrak{C},\omega_2,\ldots,\omega_n) &= \Psi_{\alpha,C,j}^{-1}\left(\mathcal{F}_{\alpha,j}(\vk\omega)\right)-\sum_{l\in{\tt S}^{\star}_j}\frac{\hat{p}_{l}\omega_l}{\hat{p}_j}\\  &= \psi_{j}^{-1}\left(\sum_{l=j+1}^{n} \left(\frac{\mathfrak{r}_j}{p_j} - \frac{\mathfrak{r}_l}{p_l}\right)p_l\omega_l\right)-\sum_{l = j}^{n}\frac{p_{l}\omega_l}{p_j},\\
    \mathcal{D}_{\alpha,j}(\mathfrak{C},\omega_2,\ldots,\omega_n) &= \Psi_{\alpha,C,j}^{-1}\left(\mathcal{F}_{\alpha,j}(\vk\omega)\right) -\sum_{l\in{\tt S}^{\star}_{j+1}}\frac{\hat{p}_{l}\omega_l}{\hat{p}_j}\\
    &= \psi_{j}^{-1}\left(\sum_{l=j+1}^{n} \left(\frac{\mathfrak{r}_j}{p_j} - \frac{\mathfrak{r}_l}{p_l}\right)p_l\omega_l\right) -\sum_{l=j+1}^{n}\frac{p_{l}\omega_l}{p_j}.
\end{align*}
Hence the claim follows due to \netheo{main_light}.\QED

\subsection{Proof related to the example of a tandem network}
\label{A5} We now set out to prove \eqref{tandem}.
In this case, for any $k\in\sprod{m}$,
\begin{align*}
    &{\tt I}_k = \sprod{q_{k+1}-1}^{q_{k}+},\quad \mathcal{r}_k(u) = r_{q_k}(u),
\end{align*}
for any $j\in\sprod{n}$,
\begin{align*}
    &\hat{p}_j = 1,\quad{\tt S}^{\star}_j = {j},
\end{align*}
and, for any $j\in\sprod{q_{k+1}-2}^{q_k+}$,
\begin{align*}
    &{\tt D}^{\star}_j = \{j+1\},\qquad {\tt D}^{\star}_{q_{k+1}-1} = \varnothing.
\end{align*}

To apply \netheo{main_light} we first calculate the constants $\mathcal{A}_{\alpha,k}(\mathfrak{C},\vk\omega)$ defined in \eqref{A_cal}: for any $k\in\sprod{m}$ and $\vk\omega\in R_+^{{\tt I}_k}$,
\begin{align*}
    \mathcal{A}_{\alpha,k}(\mathfrak{C},\vk\omega) &= \sum\limits_{l\in{\tt I}_{k}}\frac{\mathfrak{r}_{l}}{\hat{p}_{l}}\sum\limits_{m\in{\tt D}^{\star}_l}\hat{p}_m\omega_m - \sum\limits_{l\in{\tt I}_{k}}\omega_l\mathfrak{r}_l - \mathfrak{C} \left(\sum\limits_{l\in{\tt S}^{\star}_{q_k}}\hat{p}_{l}\omega_l\right)^{\alpha}\\
    &=  \sum\limits_{l=q_k}^{q_{k+1}-2}\left(\frac{\mathfrak{r}_{l}\hat{p}_{l+1}}{\hat{p}_{l}} - \mathfrak{r}_{l+1}\right)\omega_{l+1} - \omega_{q_k}\mathfrak{r}_{q_k} - \mathfrak{C} \left(\hat{p}_{q_k}\omega_{q_k}\right)^{\alpha}\\
    &=  \sum\limits_{l=q_k+1}^{q_{k+1}-1}\bigl(\mathfrak{r}_{l-1} - \mathfrak{r}_{l}\bigr)\omega_{l} - \omega_{q_k} - \mathfrak{C} \omega_{q_k}^{\alpha}.
\end{align*}

Then we calculate the constants $\mathcal{F}_{\alpha,j}(\vk\omega)$ defined in \eqref{F_cal}: for any $j\in\sprod{q_{k+1}-2}^{q_k+}$, $k\in\sprod{m}$,
\begin{align*}
    \mathcal{F}_{\alpha,j}(\vk\omega) &=\sum_{\substack{l=j+1 \\ l\in{\tt I}_{\mathcal{k}(j+1)}}}^{n}\left(\frac{\mathfrak{r}_{l-1}}{\hat{p}_{l-1}} - \frac{\mathfrak{r}_l}{\hat{p}_{l}}\right)\sum_{i\in{\tt S}^{\star}_{l}}\hat{p}_{i}\omega_i=\sum_{l=j+1}^{q_{k+1}-1}\left(\mathfrak{r}_{l-1} - \mathfrak{r}_l\right)\omega_l.
\end{align*}
Using that for any $j\in\sprod{n-1}$, according to \eqref{Psi_def},
\begin{align*}
    \Psi_{\alpha,\mathfrak{C},j}(s) = \mathfrak{r}_js + \mathfrak{C}\hat{p}_j^{\alpha}s^{\alpha} = \mathfrak{r}_js + \mathfrak{C}s^{\alpha} = \psi_j(s),
\end{align*}
we have due to \eqref{C_cal} and \eqref{D_cal} that, for any $j\in\sprod{q_{k+1}-2}^{q_k+}$, $k\in\sprod{m}$, 
\begin{align*}
    \mathcal{C}_{\alpha,j}(\mathfrak{C},\omega_2,\ldots,\omega_n) &= \Psi_{\alpha,C,j}^{-1}\bigl(\mathcal{F}_{\alpha,j}(\vk\omega)\bigr)-\sum_{l\in{\tt S}^{\star}_j}\frac{\hat{p}_{l}\omega_l}{\hat{p}_j}= \psi_{j}^{-1}\left(\sum_{l=j+1}^{q_{k+1}-1}\left(\mathfrak{r}_{l-1} - \mathfrak{r}_l\right)\omega_l\right)-\omega_j,\\
    \mathcal{D}_{\alpha,j}(\mathfrak{C},\omega_2,\ldots,\omega_n) &= \Psi_{\alpha,C,j}^{-1}\bigl(\mathcal{F}_{\alpha,j}(\vk\omega)\bigr) -\sum_{l\in{\tt S}^{\star}_{j+1}}\frac{\hat{p}_{l}\omega_l}{\hat{p}_j}= \psi_{j}^{-1}\left(\sum_{l=j+1}^{q_{k+1}-1}\left(\mathfrak{r}_{l-1} - \mathfrak{r}_l\right)\omega_l\right) -\omega_{j+1}.
\end{align*}
Hence the claim follows due to \netheo{main_light}.

\bibliographystyle{IEEEtranSN}
{\small \bibliography{EEEA}}

\end{document}